\documentclass[a4paper]{article}

\usepackage[%
  paperwidth=192mm,
  paperheight=262mm,
  vmargin={19mm,19mm},
  hmargin={13.7mm,13.7mm},
  headsep=12pt,
  footskip=12pt,
]{geometry}
\usepackage[authoryear,longnamesfirst]{natbib}
\usepackage{hyperref}
\usepackage{subcaption}
\usepackage{graphicx}
\usepackage{fancyhdr}
\usepackage{amsmath}
\usepackage{url}
\usepackage{amssymb}
\usepackage{lineno}
\usepackage{setspace}
\usepackage{cleveref}
\usepackage[inline]{enumitem}
\usepackage{amsmath}
\usepackage{mathtools, cuted}
\usepackage[ruled,linesnumbered]{algorithm2e} % For algorithms
\usepackage[T1]{fontenc}
\begin{document}
\onehalfspacing % Set line spacing to 1.5

% \let\WriteBookmarks\relax
% \def\floatpagepagefraction{1}
% \def\textpagefraction{.001}

% \shorttitle{Differentiable Smoothed Particle Hydrodynamics for Adjoint Optimization and Machine Learning}    

% \shortauthors{R. Winchenbach \& N. Thuerey}  

\title{Solving Boundary Handling Analytically in Two Dimensions for Smoothed Particle Hydrodynamics}  




% \author[1]{Rene Winchenbach}

% \cormark[1]


% \ead{rene.winchenbach@tum.de}


% \credit{Conceptualization, Methodology, Software, Validation, Writing - Original Draft, Visualization}

\author{Rene Winchenbach\\Technical University Munich \and Andreas Kolb\\University of Siegen}

% % Address/affiliation
% \affiliation[1]{organization={Technical University Munich},
%             addressline={Boltzmanstraße 3}, 
%             city={Garching bei München},
% %          citysep={}, % Uncomment if no comma needed between city and postcode
%             postcode={85748}, 
%             state={Bavaria},
%             country={Germany}}

% \author[1]{Nils Thuerey}[orcid=0000-0001-6647-8910]%[]

% % Footnote of the second author
% % \fnmark[2]

% % Email id of the second author
% \ead{nils.thuerey@tum.de}

% % URL of the second author
% % \ead[url]{}

% % Credit authorship
% \credit{Writing - Review \& Editing, Supervision, Funding acquisition}

% % Corresponding author text
% \cortext[1]{Corresponding author}


% \nonumnote{}

\maketitle
\begin{abstract}
We present a fully analytic approach for evaluating boundary integrals in two dimensions for Smoothed Particle Hydrodynamics (SPH).
Conventional methods often rely on boundary particles or wall re-normalization approaches derived from applying the divergence theorem, whereas our method directly evaluates the area integrals for SPH kernels and gradients over triangular boundaries.
This direct integration strategy inherently accommodates higher-order boundary conditions, such as piecewise cubic fields defined via Finite Element stencils, enabling analytic and flexible coupling with mesh-based solvers.
At the core of our approach is a general solution for compact polynomials of arbitrary degree over triangles by decomposing the boundary elements into elementary integrals that can be solved with closed-form solutions.
We provide a complete, closed-form solution for these generalized integrals, derived by relating the angular components to Chebyshev polynomials and solving the resulting radial integral via a numerically stable evaluation of the Gaussian hypergeometric function $_2F_1$.
Our solution is robust and adaptable and works regardless of triangle geometries and kernel functions.
We validate the accuracy against high-precision numerical quadrature rules, as well as in problems with known exact solutions.
We provide an open-source implementation of our general solution using differentiable programming to facilitate the adoption of our approach to SPH and other contexts that require analytic integration over polygonal domains.
Our analytic solution outperforms existing numerical quadrature rules for this problem by up to five orders of magnitude, for integrals and their gradients, while providing a flexible framework to couple arbitrary triangular meshes analytically to Lagrangian schemes, building a strong foundation for addressing several grand challenges in SPH and beyond.
\end{abstract}

% \includegraphics[width=\linewidth]{flowchart.drawio.png}

% \begin{highlights}
% \item A novel SPH framework built around differentiability and for Machine Learning
% \item The applicability of our framework to solve adjoint problems, e.g., solving inverse problems
% \item A novel approach to address particle-shifting, a central issue in SPH, using differentiable SPH operators
% \end{highlights}

% Keywords
% \begin{keywords}
%  \sep Differentiability\sep Adjoint Optimization\sep Machine Learning \sep Smoothed Particle Hydrodynamics \sep Solver Framework
% \end{keywords}

% Main text

\section{Introduction}
\label{sec:introduction}

Smoothed Particle Hydrodynamics~(SPH) is a Lagrangian simulation method originally proposed for the simulation of stellar phenomena in astrophysical simulations~\cite{Gingold1977SPH}.
SPH has since found a wide range of applications, including in computational fluid dynamics~\cite{vacondio2021grand,antuono2021delta}, astrophysics~\cite{hopkins2013general,frontiere2017crksph}, and computer animation~\cite{SPHTutorial,DBLP:journals/tog/WinchenbachA020}, due to the inherent Lagrangian nature of SPH leading to significant benefits for free-surface flows and complex deformations.
However, the Lagrangian nature of SPH does not lend itself to incorporating rigid boundaries or to the coupling to Eulerian mesh-based schemes, as rigid objects and mesh elements are not Lagrangian in nature.
That is, they are often described using meshes, and integrating a rigid body into an SPH simulation requires either finding an accurate Lagrangian representation for the rigid body or incorporating non-Lagrangian aspects into SPH.
Resolving this impedance mismatch between the Lagrangian nature of SPH and the discrete element-based description of other schemes and boundaries is of significant importance to both the numerical accuracy and versatility of SPH in computational fluid dynamics and beyond.

Particle-based rigid boundary handling methods can be distinguished into ghost particle methods~\cite{colagrossi2003numerical,libersky1993high,marrone2011delta,yildiz2009sph,english2022modified}, which either place virtual particles inside the boundary that cover the entire support domain of a particle, or on the surface of the boundary~\cite{Akinci2012RigidFluid,band2017moving,monaghan1994simulating}, where either approach may utilize a consistent set of particles or generate them dynamically.
Particle-based methods generally suffer from sampling problems, e.g., placing particles on the surface in a uniform manner, and struggle with representing flat geometries as the particle sampling can introduce artifacts in the density field~\cite{band2017moving}. 
Consequently, corrective terms are required to counteract these problems~\cite{Akinci2012RigidFluid,band2018pressure}, e.g., pseudo-volumes accounting for deficiencies in the particle sampling or to represent the rigid body using boundary surface particles.

Instead of finding a Lagrangian representation of the boundary, integral formulations directly incorporate a non-Lagrangian boundary representation into the SPH formulation.
The most common approach for boundary integrals is determining a wall renormalization factor $\gamma$, originally introduced by Kulasegaram et~al.~\cite{kulasegaram2004variational}, where the integrals themselves are commonly evaluated using the divergence theorem.
There has been significant research in this direction by incorporating gradient operators~\cite{ferrand2013unified}, using semi-analytic~\cite{marongiu2010free, deleffe2011normal}, numerical~\cite{deleffe2011normal, chiron2019fast} or fully analytic solutions to the integrals~\cite{violeau2009exact,leroy2014unified,mayrhofer2015unified}, for both non-viscous and viscous flows in both 2D and 3D.
Instead of using wall renormalization approaches, an alternative approach is the direct evaluation of the contribution of a boundary object to a fluid particle, enabling a straightforward integration into fluid simulation models.
These methods work by either numerically evaluating the contributed density~\cite{koschier2017density} or a representative volume term~\cite{bender2019volume} or by semi-analytically evaluating the kernel function over the boundary domain~\cite{DBLP:journals/tog/WinchenbachA020}.
However, these schemes do not generally provide exact solutions and rely on approximations that may not always be applicable.

In this paper, we overcome these limitations by introducing a novel and generalized framework for the fully analytic, direct evaluation of 2D SPH boundary integrals over triangular domains.
Our solution can be applied not just to handling rigid boundaries, but is also capable of coupling SPH exactly to mesh-based simulation schemes, removing a significant gap in the fundamental difference in discretization.
Using our general solution, we can find these closed-form integrals for any compact polynomial kernel, e.g., B-Splines and Wendland kernels, of arbitrary order.
The core of our approach relies on a decomposition of triangles into fundamental forms that can be analytically computed for only a small set of terms, i.e., $r^k$, $r^k\cos n\theta$, and $r^k\cos n\theta \sin m\theta$.
We then establish the closed-form analytic solutions for these general integrals by relating the angular components to Chebyshev polynomials and solving the resulting radial integrals through a numerically stable evaluation of the Gaussian hypergeometric function $_2F_1$.
This approach inherently supports higher-order boundary conditions, such as piecewise linear fields defined by barycentric interpolation, enabling more sophisticated and accurate physical interactions.

The remainder of this paper is organized as follows.
Section~\ref{sec:SPH} briefly reviews the SPH formulation and defines our general boundary integral formulation.
Section~\ref{sec:barycentric} describes the integration of piecewise polynomial fields using barycentric interpolation into the boundary integral formulation and defines the fundamental integral terms that need to be solved.
Section~\ref{sec:integralTerms} covers the decomposition of arbitrary triangles into fundamental forms, i.e., sectors, segments, and stubs, with well-defined integration bounds and provides the general algorithm of combining these fundamental terms into integrals over actual triangles.
Section~\ref{sec:integrals} describes the fully-analytic closed-form solution of the integral terms described in Section 4 and the evaluation of the $_2F_1$ function for piecewise constant and polynomial boundary terms.
Section~\ref{sec:eval} evaluates the numerical accuracy of our approach by comparing the results to high-order numerical quadratures, which rely on the same description as in Section 3, as well as evaluating the accuracy against problems with known exact solutions.
This section also demonstrates the application of our scheme to a simple SPH simulation.
Finally, Section~\ref{sec:conclusion} provides concluding remarks and discusses future work.
Note that we provide an open-source implementation of the framework both in C++ and PyTorch, enabling the use of our solution in SPH and beyond, with the PyTorch approach enabling novel applications using differentiable, exact, higher-order coupling to boundaries and mesh-based solvers.

\section{Governing Equations and Boundary Integral Formulation}
\label{sec:SPH}
SPH fundamentally is based on an identity of a function, where a field quantity $A$ at any point $\mathbf{x}\in\mathbb{R}^2$ can be evaluated as an integral over the entire input space~\cite{price2012smoothed}:
\begin{equation}
\langle A(\mathbf{x})\rangle = \int_{\mathbb{R}^2} A(\mathbf{x}^\prime)\delta(\mathbf{x},\mathbf{x}^\prime) d\mathbf{x}^\prime,
\end{equation}
where $\delta$ is the Dirac delta function.
Note that this definition works in any number of spatial dimensions, but we only discuss two-dimensional SPH here, as our method is only applicable in 2D.
The Dirac delta function is then replaced by a compact kernel function $W$ with a support radius $h$~\cite{price2012smoothed} yielding an approximation for a field quantity $A$ at any point $\mathbf{x}$ based on a compact spherical domain $\Omega$, with cutoff radius $h$ and centered at $\mathbf{x}$, as
\begin{equation}
\label{eqn:sphintegral}
\langle A(\mathbf{x})\rangle \approx \int_\Omega A(\mathbf{x}^\prime)W(\mathbf{x}-\mathbf{x}^\prime, h) d\mathbf{x}^\prime.
\end{equation}
This can also be applied to gradients of generic field quantities to yield a naïve continuous gradient formulation as
\begin{equation}
\label{eqn:nablaIntegral}
\langle \nabla A(\mathbf{x})\rangle \approx \int_\Omega A(\mathbf{x}^\prime)\nabla_\mathbf{x}W(\mathbf{x}-\mathbf{x}^\prime, h) d\mathbf{x}^\prime.
\end{equation}

SPH then uses particles to discretize the fluid domain, which carry discrete values of the field quantities, where the apparent volume $V_j = \frac{m_j}{\rho_j}$ of a particle is used as weights~\cite{monaghan2005smoothed}, yielding
\begin{equation}
\langle \nabla A(\mathbf{x})\rangle = \sum_j \frac{m_j}{\rho_j} A_j W(\mathbf{x} - \mathbf{x}_j,h).
\end{equation}
Here, $j$ refers to the neighboring particles around $\mathbf{x}$, i.e., all particles $j$ with $|\mathbf{x} - \mathbf{x}_j| \leq h$, with $m_j$ and $\rho_j$ being the mass and density of the particle, respectively. 
Gradient terms can analogously be determined, however, such terms would neither result in symmetric forces nor be exact for constant functions~\cite{price2012smoothed} in the presence of particle disorder.
A gradient formulation that is exact for constant functions, but not symmetric, can be derived as
\begin{equation}
\begin{split}
\label{eqn:nablasIntegralDifference}
\nabla A &= \frac{1}{\rho}\left[\langle\nabla(\rho A)\rangle - A\langle\nabla\rho\rangle\right],
\end{split}
\end{equation}
also called the difference formulation. A gradient formulation that is symmetric, but not exact for constant functions, can be derived as
\begin{equation}
\begin{split}
\label{eqn:nablasIntegralSymmetric}
\nabla A &= \rho\left[\left<\nabla\left(\frac{A}{\rho}\right)\right> + \frac{A}{\rho^2}\langle\nabla\rho\rangle\right],
\end{split}
\end{equation}
also called symmetric formulation~\cite{price2012smoothed}. The discretized forms of Eq.~\eqref{eqn:nablasIntegralDifference} and Eq.~\eqref{eqn:nablasIntegralSymmetric} yield:
\begin{equation}
\label{eqn:sphdiffops}
\begin{split}
\langle\nabla A(\mathbf{x})\rangle &\approx \frac{1}{\rho(\mathbf{x})}\sum_j m_j (A_j - A(\mathbf{x}))\nabla_\mathbf{x}W(\mathbf{x} - \mathbf{x}_j,h)\\
\langle\nabla A(\mathbf{x})\rangle &\approx {\rho(\mathbf{x})}\sum_j m_j \left(\frac{A_j}{\rho_j^2} - \frac{A(\mathbf{x})}{\rho(\mathbf{x})^2}\right)\nabla_\mathbf{x}W(\mathbf{x} - \mathbf{x}_j,h).
\end{split}
\end{equation}
Note that analogous formulations exist for both divergence $\nabla\cdot A(\mathbf{x})$ and curl $\nabla\times A(\mathbf{x})$~\cite{monaghan1994simulating}.
For our boundary integral evaluation, we require the continuous forms of these terms, which can be directly defined using the continuous form of the SPH gradient operator as 
\begin{align}
\label{eqn:intsphdiffops}
\nabla A \approx &\frac{1}{\rho(\mathbf{x})}\left[\langle\nabla(\rho A)\rangle - A\langle\nabla\rho\rangle\right] 
&&=\frac{1}{\rho(\mathbf{x})}\int_\Omega \rho(\mathbf{x}^\prime)\left[A(\mathbf{x}^\prime) - A(\mathbf{x}) \right]\nabla_\mathbf{x}W(\mathbf{x}-\mathbf{x}^\prime,h) d\mathbf{x}^\prime \\
\nabla A \approx &\rho(\mathbf{x})\left[\left<\nabla\left(\frac{A}{\rho}\right)\right> + \frac{A}{\rho^2}\langle\nabla\rho\rangle\right] 
&&={\rho(\mathbf{x})}\int_\Omega \rho(\mathbf{x}^\prime) \left[\frac{A(\mathbf{x}^\prime)}{\rho^2(\mathbf{x}^\prime)} - \frac{A(\mathbf{x})}{\rho^2(\mathbf{x})}\right]\nabla_\mathbf{x}W(\mathbf{x}-\mathbf{x}^\prime,h) d\mathbf{x}^\prime .
\end{align}

A general definition of the kernel function, based on Dehnen and Aly~\cite{dehnen2012improving}, is defined as 
\begin{equation}
    W(\mathbf{x},h) = \frac{1}{h^d} C_d k\left(\frac{|x|}{h}\right),
\end{equation}
with $k(q)$ being the actual kernel, $d$ the dimension and $C_d$ a normalization factor.
Note that, without loss of generality for the derivation, we can assume the smoothing scale and the cut-off radius to be identical.
Consequently, the derivative of the kernel function is given as
\begin{equation}
    \nabla_\mathbf{x} W(\mathbf{x}- \mathbf{x}^\prime,h) = \frac{\mathbf{x}-\mathbf{x}^\prime}{|\mathbf{x}-\mathbf{x}^\prime|}\frac{1}{h^{d+1}} C_d \frac{\partial k\left(\frac{|\mathbf{x}-\mathbf{x}^\prime|}{h}\right)}{\partial\frac{|\mathbf{x}-\mathbf{x}^\prime|}{h}}.
\end{equation}
We use the Wendland 4 kernel with $C_2=\frac{27}{16\pi}$ and $k(q) = \left[1-q\right]^6_+ \left(1+6q+\frac{35}{3}q^2\right)$, with $\left[\cdot\right]_+=\max\left\{\cdot,0\right\}$. 
As the kernel can be solely expressed using integer powers of the radial distance $q$ within the support domain $\Omega$, i.e., it can be written as $k(q) = \sum_i b_k q^k$, solving the boundary integral formulation for an arbitrary power $q^k$ yields a solution to the entire kernel function.
This can be extended to piecewise kernel functions trivially, such as the B-Splines, through careful evaluation of the integral bounds for each piece of the kernel function.

We now assume that there is a boundary domain $\Omega_b$, consisting of individual triangle elements $\mathcal{T}$, with no requirements placed on the individual triangles' shapes and connectivity.
Note that our approach still requires disjoint triangles, i.e., triangle meshes with no overlap, to yield physically plausible answers.
Given a set of neighboring triangles $\mathcal{N}_\mathcal{T}$, i.e., triangles that overlap the support radius $h$ relative to a point in space $\mathbf{x}$, the integral over the entire boundary domain is evaluated as
\begin{equation}
\label{eqn:boundaryint}
    \langle A(\mathbf{x})\rangle \approx \sum_{\mathcal{T}\in\mathcal{N}_\mathcal{T}}\int_\mathcal{T}A(\mathbf{x}^\prime)W(\mathbf{x}-\mathbf{x}^\prime,h)d\mathbf{x}^\prime.
\end{equation}
As we assume triangles to be disjointed, we can treat each triangle that is to be integerated independently as their connectivity to neighboring triangles does not affect our scheme.
To further simplify the integral solutions, we first transform this equation into a normalized form, i.e., we remove unnecessary degrees of freedom from the equation, without loss of generality.
This involves transforming the coordinate system such that the cutoff radius is equal to $1$, note that this does not necessarily mean that the smoothing length is $1$, and that the position $\mathbf{x}$ at which this integral is evaluated is at the origin $\mathbf{0}$. 
%
%Furthermore, we also require an inclusion of barycentric interpolation and of difference and symmetric gradient formulations.
%
Each individual triangle $\mathcal{T}$ consisting of three vertices located at positions $\mathbf{\nu}_1,\mathbf{\nu}_2$ and $\mathbf{\nu}_3\in\mathbb{R}^2$ is thus shifted by $\mathbf{x}$ and inversely scaled by $h$ to ensure that $h=1$ for the purposes of integration by applying a homogeneous transformation matrix 
$$
C = \left[\begin{array}{ccc}
~1/h&~0&~-\frac{\mathbf{x}_x}{h}\\
0&~1/h&~-\frac{\mathbf{x}_y}{h}\\
0&~0&1
\end{array}\right],
$$ 
which yields a transformed triangle $\mathcal{T}^\star$, and vertices $\nu_i^\star = C \nu_i$, with
\begin{equation}\begin{split}
    \langle A(\mathbf{0})\rangle &\approx 
    \int_{\mathcal{T}^\star}A(\mathbf{x}^\prime)W(-\mathbf{x}^\prime,1)d\mathbf{x}^\prime,\\
    \langle\nabla A(\mathbf{0})\rangle &\approx 
    \frac{1}{h}\int_{\mathcal{T}^\star}A(\mathbf{x}^\prime)\nabla_\mathbf{x}W(-\mathbf{x}^\prime,1)d\mathbf{x}^\prime,
    \end{split}
\end{equation}
Note that this does not change the result of the kernel integral as the kernel function $W$ is invariant to translation, i.e., only relative positions matter. 
For the gradient term, the additional division by the support radius $h$ is required as part of the u-substitution.
Referring back to the general kernel formulation~\cite{dehnen2012improving}, this means that the kernel function $W(\mathbf{x},1)$ can be replace with $C_d k(|\mathbf{x}|)$, where we drop the $C_d$ as this term can simply be multiplied back onto the result after integration, without loss of generality. 
This yields an integral for the boundary contribution as
\begin{equation}
\label{eqn:singleTriangle}
    \langle A(\mathbf{0})\rangle \approx \sum_{\mathcal{T}^ \star\in\mathcal{N}_\mathcal{T}}\int_{\mathcal{T}^ \star}A(\mathbf{x}^\prime)w(|-\mathbf{x}^\prime|)d\mathbf{x}^\prime,
\end{equation}
with $w(|\mathbf{x}|) = W(\mathbf{x},1)$.
To simplify the integration process later on, we apply an improper rotation to a triangle to transform a given triangle into a triangle for which a solution is known.
While Eq.~\ref{eqn:singleTriangle} is invariant to rotations of the coordinate system, integrals of gradient terms need to be adjusted when an improper rotation matrix $R$ is applied to the local coordinate system, after scaling and shifting, to yield a new triangle $\mathcal{T}^R$, and vertices $\mathbf{\nu}^r_i = R \mathbf{\nu}^\star_i$. This, after dropping the subscript of $\nabla$ for readability, yields a term that needs to be transformed with the inverse of $R$, i.e., \begin{equation}
\label{eqn:improper}
    \int_{\mathcal{T}^\star}A(\mathbf{x}^\prime)\nabla w(|-\mathbf{x}^\prime|)d\mathbf{x}^\prime = R^{-1} \int_{\mathcal{T}^R}A(\mathbf{x}^\prime)\nabla w(|-\mathbf{x}^\prime|)d\mathbf{x}^\prime.
\end{equation}

To evaluate the integral for difference gradient formulations (Eq.~\ref{eqn:intsphdiffops}), we first evaluate quantities $f_i = \rho_i (A_i - A(\mathbf{0}))$ for each vertex of a triangle and then perform a barycentric interpolation of these quantities across the triangle during integration.
Consequently, the overall integrals remain unchanged and only require multiplication by a different factor.
Symmetric gradient formulations can be evaluated analogously by evaluating $f_i = \rho_i\left[\frac{A_i}{\rho_i^2}+\frac{A(\mathbf{0})}{\rho(\mathbf{0})^2}\right]$ prior to integration.
Finally, divergence terms and curl terms for vector-valued fields becomes
\begin{equation}
\begin{split}
\langle \nabla\cdot A(\mathbf{x})\rangle &= \langle \nabla A_x(\mathbf{x})\rangle_x + \langle \nabla A_y(\mathbf{x})\rangle_y,\\
\langle \nabla\times A(\mathbf{x})\rangle &= \langle \nabla A_y(\mathbf{x})\rangle_x - \langle \nabla A_x(\mathbf{x})\rangle_y.
\end{split}
\end{equation}

\section{Piecewise Polynomial Field Quantities}
\label{sec:barycentric}
In Eulerian simulations using triangular meshes the two most common forms of quantities are either per-element information, i.e., piece-wise constant, or per-vertex information, i.e., piece-wise linear information.
Note that in many prior works, a common choice in SPH is assuming $A$ to be constant over each boundary element, i.e., $A(\mathbf{x}) = \text{const}:\forall\mathbf{x}\in\mathcal{T}$, i.e., limiting A to be piecewise constant.
Alternatively, using barycentric interpolation to evaluate $A$ across a triangle $\mathcal{T}$, the fields defined over boundaries can be piecewise polynomial using Finite Element Method stencils.
Replacing the generic field $A(\mathbf{x}$) with a barycentric evaluation $A_\mathcal{T}(\mathbf{x})$, over the triangle $\mathcal{T}$, and by transforming from cartesian to polar coordinates yields
\begin{equation}
\label{eqn:singleSpherical}
    \begin{split}
        \langle A(\mathbf{0})\rangle \approx \int_{\mathcal{T}}  A_{\mathcal{T}}\left(r \cos\theta, r\sin\theta\right) w(r)r d\theta dr.
    \end{split}
\end{equation}
$A_{\mathcal{T}}$ is based on a barycentric interpolation using the vertex positions $\mathbf{\nu}_i=[x_i,y_i]^T$ and their representative field values $A_i$. 
Barycentric interpolation can be defined using the canonical interpolation weights
\begin{equation}
\lambda_0(x,y) = \frac{(y_1 -y_2)(x-x_2) + (x_2 - x_1)(y-y_2)}{(y_1-y_2)(x_0-x_2)+(x_2-x_1)(y_0-y_2)},\quad
\lambda_1(x,y) = \frac{(y_2 -y_0)(x-x_2) + (x_0 - x_2)(y-y_2)}{(y_1-y_2)(x_0-x_2)+(x_2-x_1)(y_0-y_2)},
\end{equation}
with the third weight being $\lambda_2(x,y) = 1 - \lambda_0 - \lambda_1$.
This then yields the canonical barycentric interpolation for a triangle as 
\begin{equation}
A_{\mathcal{T}}(x,y)= \lambda_0(x,y)A_0 + \lambda_1(x,y)A_1 + \lambda_2(x,y)A_2 =
A_2 + \lambda_0(x,y)\left[A_0 - A_2\right] + \lambda_1(x,y)\left[A_1 - A_2\right],
\end{equation}
Note that this formulation is correct even outside the triangle and can also be used to perform point in triangle checks that will later be required for some steps of the algorithm to compute an integral.
Inserting this back into~\cref{eqn:singleSpherical} yields
\begin{equation}\langle A(\mathbf{0})\rangle \approx 
A_2\int_\mathcal{T}w(r)r d\theta dr +\left[A_0-A_2\right]\int_\mathcal{T}\lambda_0(x,y)w(r)r d\theta dr
+\left[A_1-A_2\right]\int_\mathcal{T}\lambda_1(x,y)w(r)r d\theta dr,
\end{equation}
where the first term is identical to the term that arises when assuming that the field $A$ is constant over the boundary element, i.e., $A_0=A_1=A_2=\text{const}$. 
Back inserting the definition of the barycentric interpolation with $D = {(y_1-y_2)(x_0-x_2)+(x_2-x_1)(y_0-y_2)}$ and refactoring common terms then yields
\begin{equation}
\langle A(\mathbf{0})\rangle \approx 
A_2\int_\mathcal{T}w(r)r d\theta dr
+\left[\Lambda^A_0+\Lambda^A_2\right] \int_\mathcal{T} \left[r\cos\theta - x_2\right]w(r)r d\theta dr +\left[\Lambda^A_1+\Lambda^A_3\right] \int_\mathcal{T} \left[r\sin\theta - y_2\right]w(r)r d\theta dr,
\end{equation}
with a new set of factors, which are independent of the coordinates, as
\begin{equation}
    \Lambda^A_0 = \frac{\left(y_1-y_2\right)\left(A_0-A_2\right)}{D}, \Lambda^A_1 = \frac{\left(x_2-x_1\right)\left(A_0-A_2\right)}{D},
    \Lambda^A_2 = \frac{\left(y_2-y_0\right)\left(A_1-A_2\right)}{D}, \Lambda^A_3 = \frac{\left(x_0-x_2\right)\left(A_1-A_2\right)}{D}.
\end{equation}
This can be further refactored into
\begin{equation}
\label{eqn:taus}
\begin{split}
\langle A(\mathbf{0})\rangle \approx  &
\underbrace{\left[\Lambda^A_0+\Lambda^A_2\right]}_{\tau^A_1}\int_\mathcal{T}w(r)r^2\cos\theta d\theta dr
+\underbrace{\left[\Lambda^A_1+\Lambda^A_3\right]}_{\tau^A_2}\int_\mathcal{T}w(r)r^2\sin\theta d\theta dr\\
&+\underbrace{\left[A_2 - \left(\Lambda^A_0 x_2 + \Lambda^A_1 y_2 + \Lambda^A_2 x_2 + \Lambda^A_3 y_2\right)\right]}_{\tau^A_3}\int_\mathcal{T}w(r)r d\theta dr.
\end{split}
\end{equation}
Using the earlier observation that we can replace a kernel for polynomial kernel functions with a summation over polynomials and their corresponding factors, we can then write
\begin{equation}
\langle A(\mathbf{0})\rangle \approx\sum_{k=0} b_k \left(\tau_1^A\int_\mathcal{T}r^{k+2}\cos\theta d\theta dr + \tau_2^A\int_\mathcal{T}r^{k+2}\sin\theta d\theta dr + \tau_3^A \int_\mathcal{T}r^{k+1}d\theta dr\right).
\end{equation}
Note that for this term only the third integral is required for piecewise-constant boundary integrals with $\tau_3^A=A$.
To determine gradient terms, we first consider a basic gradient formulation, where after inserting the kernel derivative, we get for a single triangle
\begin{equation}
\label{eqn:basicdiff}
\langle\nabla A\rangle \approx-\int_\mathcal{T} A_\mathcal{T}(r\cos\theta,r\sin\theta)\left[\begin{array}{c}~\cos\theta\\\sin\theta
\end{array}\right]\frac{\partial w(r)}{\partial r} r d\theta dr,
\end{equation}
Note that in polar coordinates the gradient, i.e., $\nabla|\mathbf{x}|$, becomes the $\cos$ and $\sin$ of the angle.
Using a derivation for the components of the gradient, similar to Eq.~\ref{eqn:taus}, we get equivalent integral terms by simply replacing $w(r)$ with $\cos\theta\frac{\partial w(r)}{\partial{r}}$ and $\sin\theta\frac{\partial w(r)}{\partial{r}}$. Accordingly, we find the following component-wise integrals for the basic gradient operator~\cref{eqn:basicdiff} as
\begin{equation}
    \begin{split}
&\langle\nabla A\rangle_x \approx\sum_{k=0} \frac{k}{h} b_k \left(\tau^A_1\int_\mathcal{T} r^{k+1}\cos^2\theta d\theta dr+\tau^A_2 \int_\mathcal{T}r^{k+1}\cos\theta\sin\theta d\theta dr+\tau^A_3\int_\mathcal{T} r^{k}\cos\theta d\theta dr\right),
    \end{split}
\end{equation}
and 
\begin{equation}
    \begin{split}
&\langle\nabla A\rangle_y \approx \sum_{k=0} \frac{k}{h} b_k \left(\tau^A_1\int_\mathcal{T} r^{k+1}\cos\theta\sin\theta d\theta dr+\tau^A_2 \int_\mathcal{T}r^{k+1}\sin^2\theta d\theta dr+\tau^A_3\int_\mathcal{T} r^{k}\sin\theta d\theta dr\right).
    \end{split}
\end{equation}
Note that using canonical trigonometric identities we can replace $\cos^2\theta$ with $\frac{1}{2}\left(1+\cos2\theta\right)$ and $\sin^2\theta$ with $\frac{1}{2}\left(1-\cos2\theta\right)$, which means that in total there are nine integrals, see Table~\ref{tab:integralterms}. 
Considering this set of nine integrals, there are only 4 unique integrals:
\begin{align}
g_k&= \int_\mathcal{T} x^k d\theta dx,\;&&k\geq 0,\\
g_k^{n,0} &= \int_\mathcal{T}x^k\cos\left(n\theta\right) d\theta dx,\;&&k\geq 0,n\geq1,\\
g_k^{0,m} &= \int_\mathcal{T}x^k\sin\left(m\theta\right) d\theta dx,\;&&k\geq 0,m\geq1,\\
g_k^{n,m} &= \int_\mathcal{T}x^k\cos\left(n\theta\right)\sin\left(m\theta\right) d\theta dx,\;&&k\geq 0,n\geq1,m\geq1.
\end{align}
Note that in a general context, $g_k^{n,0}$ for $k=0$ would equal $g_k$, however, in practice, this leads to issues in automated integration, and they require significantly different treatment, so we chose to keep both terms.
However, in all cases we investigated, $g_k^{0,m}$ can be treated identically to the more general $g_k^{n,m}$, reducing the number of integrals that need to be solved to $3$.
It is also important to note, that, if evaluated naïvely, integrals like $g_k^{n,m}$ can lead to inconsistencies if $n=m$; here, special care needs to be taken.
Finally, this process can be done for arbitrary derivatives of the integrals, which, using trigonometric identities, will always break down to the same set of three basic integral terms.
Furthermore, this process can be applied to any polynomial that is defined in barycentric coordinates, which includes many widely utilized Finite Element Method~(FEM) stencils.
Consequently, this means that our approach can be trivially expanded from piecewise-linear to piecewise-n boundary integrals using FEM approaches and, furthermore, this directly implies that our solution allows for the exact (both analytic and algebraic) coupling of an FEM method to SPH.
The challenge now becomes integrating these three unique terms over arbitrary triangles, which will be done in the next Section.

\begin{table*}[t]
\centering
\begin{tabular}{r|c|c|c}
                              & I & II & III \\\hline
$\cdot$                       & $\int_\mathcal{T}r^{k+2}\cos\theta d\theta dr$  & $\int_\mathcal{T}r^{k+2}\sin\theta d\theta dr$   &  $\int_\mathcal{T}r^{k+1} d\theta dr$   \\\hline
$\frac{\partial}{\partial x}$ & 
$\frac{1}{2}\left(\int_\mathcal{T}r^{k+1}d\theta dr + \int_\mathcal{T}r^{k+1}\cos2\theta d\theta dr\right)$  & 
$\int_\mathcal{T}r^{k+1}\cos\theta\sin\theta d\theta dr$   &  
$\int_\mathcal{T} r^{k}\cos\theta d\theta dr$     \\\hline
$\frac{\partial}{\partial y}$ & 
$\frac{1}{2}\left(\int_\mathcal{T}r^{k+1}d\theta dr - \int_\mathcal{T}r^{k+1}\cos2\theta d\theta dr\right)$  & 
$\int_\mathcal{T}r^{k+1}\sin^2\theta d\theta dr$   &  
$\int_\mathcal{T} r^{k}\sin\theta d\theta dr$    
\end{tabular}
\caption{Integral terms for normal and gradient terms, for piecewise-constant (III) and piecewise-polynomial (I-III) interpolations}
\label{tab:integralterms}
\end{table*}

\section{Integral Terms}
\label{sec:integralTerms}

In the previous section, we identified 3 integral terms which need to be evaluated for each boundary triangle. 
However, directly evaluating these terms over an arbitrary triangle is not feasible.
Instead, we propose a decomposition approach that combines integrals for geometries which can be integrated directly.
The basic geometries in our case are \begin{enumerate*}[label=\upshape(\itshape\alph*\upshape)]
    \item cones,
    \item segments and
    \item triangles with one vertex at the origin and at least one other vertex in the support domain, referred to as stubs
\end{enumerate*}. 
These basic geometries are then combined to create 
\begin{enumerate*}[label=\upshape(\itshape\alph*\upshape)]
    \item arbitrary wedges, i.e., triangles with none or one vertex inside the support domain,
    \item arbitrary triangles with two vertices in the support domain and
    \item arbitrary triangles with three vertices in the support domain
\end{enumerate*}. 
To simplify the integrals, it is helpful to only treat a single spatial configuration, e.g., circular regions where one edge of the triangle coincides with the x-axis.
In this regard, the $g_k^{0,0}$ term can be treated straightforwardly, as it is radially symmetric, whereas all other terms require an improper rotation first to transform the spatial configuration to a single base case; see~\cref{eqn:improper}.
Consequently, we can focus our derivations on a single spatial configuration, as long as any arbitrary shape of the same kind can be transformed into this configuration using an improper rotation. 
For brevity, we will discuss the evaluation of the integrals for an arbitrary function $g(r,\theta)$, i.e., any choice from the functions in Tab.~\ref{tab:integralterms} as the important aspect in this section is finding the right integration bounds. Later on, Section~\ref{sec:integrals} will cover evaluating the integrals over the given bounds and numerical challenges in evaluating the analytic forms.
Accordingly, we will consider an integral of the form
\begin{equation}
    \begin{split}
        \int_\mathcal{T} g(r,\theta) d\theta dr.
    \end{split}
\end{equation}
Note that this formulation is independent of the barycentric interpolation, as these are just constant pre-factors of the integrals, and independent of the used kernel function, as we treat $g$ as a black-box function for the purposes of this section.
However, a very important caveat to this is that the weights of the barycentric interpolation and the results of gradient terms need to be potentially rotated using a given improper rotation matrix.
Accordingly, the processes developed here need to perform the corresponding transformations, even though the actual integrals are not.

\subsection{Integrals Over Discs}
\label{sec:triangles:cones}

The first set of integral bounds we need to consider is fundamental integrals over circular sectors, i.e., cases where one edge of the triangle is completely outside the circle and one vertex is at the origin.
This case occurs as a subproblem in some other cases and in its most general form requires the integration over an opening from angle $\alpha$ to $\beta$ and over a limited radius $d\leq1$.
There are two noteworthy special cases of this for integrating over the entire support domain, i.e., from $0$ to $2\pi$, and integrals over half the domain, i.e., $-\pi/2$ to $\pi/2$.
The former case never requires an adjustment as the result is correct even for all potential geometries, whereas the case of a half support integral and the more general integral between two angles requires an adjustment of the input such that the angle is symmetric with respect to the x-axis.
This means that our first integral of concern is 
\begin{equation}
    \begin{split}
        \iint_\mathcal{C} g(r,\theta) r d\theta dr = \int_0^d\int_0^{2\pi}g(r,\theta) d\theta dr.
    \end{split}
\end{equation}
which leads to a straightforward evaluation, see Algorithm~\ref{algo:circle}.
The second case, an integral over half the support domain, simply requires changing the integration bounds to yield
\begin{equation}
    \begin{split}
        \iint_\mathcal{\text{Half}} g(r,\theta) r d\theta dr = \int_0^d\int_{-\frac{\pi}{2}}^{\frac{\pi}{2}}g(r,\theta) d\theta dr.
    \end{split}
\end{equation}
Note that this integral is only required as part of a later integral term and does not require an explicit evaluation process.
The final variant of a general sector, the integral covers the range $[\gamma^0, \gamma^0 + \gamma]$ and the radius $[0,\min\{d,1\}]$. 
To simplify the integral terms, we assume that the lower angle $\gamma^0$ is equal to $0$, which can be ensured for any sector through rotation by $-\gamma^0$. Consequently, we can find a function $\mathcal{S}^\mathcal{T}$ dependent on an opening angle of the cone $\gamma$ and a radial limit $d$ as
\begin{equation}
\label{eqn:conical}K^\text{sector}(\gamma,d) = 
    \mathcal{S}^\mathcal{T}(\gamma,d) = \int_0^d\int_0^\gamma g(r,\theta)d\theta dr.
\end{equation}

\DontPrintSemicolon
\begin{algorithm}[t]
    \SetAlgoLined
    \textbf{Function} sphericalIntegral(d, $\mathcal{T}$, $f$):\;
    \Indp
        Calculate $\tau^f_1,\tau^f_2,\tau^f_3$ for $\mathcal{T}$; see~\cref{eqn:taus}\;
        \textbf{Return} $\mathcal{C}^\mathcal{T}(d), \left[\mathcal{C}^\mathcal{T}_x(d),\mathcal{C}^\mathcal{T}_y(d)\right]^T$\;
    \Indm
\caption{Procedure to evaluate integrals over discs. This procedure requires as input the radius $d$ of the disc, the reference triangle $\mathcal{T}$, and the scalar field $f$. The procedure returns the integral over the region and the gradient.}
\label{algo:circle}
\end{algorithm}

For a given triangle $\mathcal{T}$, whose intersection with the support domain results in a sector, we first determine the orientation of $\mathbf{\nu}_1$, i.e., $\alpha = \operatorname{atan2}(y_1,x_1)$ and the sign of the y-coordinate of rotating $\mathbf{\nu}_2$ by $\alpha$ as $s=\operatorname{sgn}(x_2 \sin-\alpha + y_2 \cos-\alpha)$ to construct an improper rotation matrix $C$ with an angle $-\alpha$ and a mirroring of the y-axis afterwards if $s<0$. 
Accordingly, in in the transformed coordinate system the integral bounds go angle-wise from $0$ to $\gamma$ with $\gamma=\operatorname{atan2}(y^\prime_2,x^\prime_2)$, for $\mathbf{\nu}_2^\prime$ being the result of applying $C$ to $\mathbf{\nu}_2$, and radially from $0$ to $d$. Note that the reference triangle $\mathcal{T}$ also needs to be transformed with $C$. 
Algorithm~\ref{algo:sector} shows how this process can be realized.
To actually solve the integrals, we can simply evaluate the most general form, i.e., 
\begin{equation}
    \iint_\text{Sector}g(k,\theta)d\theta dr = \int_0^d\int_\alpha^\beta g(k,\theta)d\theta dr.
\end{equation}

\DontPrintSemicolon
\begin{algorithm}[t]
    \SetAlgoLined
    \textbf{Function} sectorIntegral($\mathbf{\nu}_1$, $\mathbf{\nu}_2$, d, $\mathcal{T}$, $f$):\;
    \Indp
        $\alpha_1 \leftarrow \operatorname{atan2}(y_1, x_1)$\;
        $s \leftarrow \operatorname{sgn}\left(x_2 \sin \alpha_1 + y_2 \cos -\alpha_1\right)$ \;
        $ C \leftarrow \left[\begin{array}{ccc}
~\cos\left(-\alpha_1\right) & \sin\left(\alpha_1\right) & 0\\
-s\sin\left(\alpha_1\right) & s\cos\left(\alpha_1\right) & 0\\
0 &0 & 1
\end{array}\right],$\;
$\mathcal{T}^\prime \leftarrow C \cdot \mathcal{T}, \mathbf{\nu}_1^\prime \leftarrow C \cdot \mathbf{\nu}_1 , \mathbf{\nu}_2^\prime \leftarrow C \cdot \mathbf{\nu}_2 $\;
$\gamma \leftarrow \operatorname{atan2}(\nu^\prime_{2,y}, \nu^\prime_{2,x})$\;
        Calculate $\tau^f_1,\tau^f_2,\tau^f_3$ for $\mathcal{T}^\prime$; see~\cref{eqn:taus}\;
        \textbf{Return} $S^{\mathcal{T}^\prime}(d,\gamma), T^{-1}\cdot\left[S^{\mathcal{T}^\prime}_x(d,\gamma),S^{\mathcal{T}^\prime}_y(d,\gamma)\right]^T$\;
    \Indm
\caption{Procedure to evaluate integrals over sectors. This procedure requires as an input the radius of the sector region $d$, the reference triangle $\mathcal{T}$, the scalar field $f$ and two vertices $\mathbf{\nu}_1$ and $\mathbf{\nu}_2$ that describe the sector, note that $\mathbf{\nu}_0$ is implied to be $\mathbf{0}$. The procedure returns the integral over the region and the gradient.}
\label{algo:sector}
\end{algorithm}

\subsection{Integrals Over Segments}
\label{sec:triangles:segments}

For circle segments, which are the result of one edge, $\overline{\mathbf{\nu}_1\mathbf{\nu}_2}$ by convention, intersecting the support domain and no other edge intersecting the domain, the integral depends on the signed distance of the segment, from the origin, and the relative orientation of the segment. 
We simplify this configuration by assuming that the segment is orthogonal to the x-axis and that its x-coordinates are equal to the signed distance.
We determine the opening angle of the segment with distance $r$, $d\leq r \leq 1$, as $\gamma(r) = 2 \arccos\frac{r}{d}$ using canonical trigonometry.
The angular limits of integration for a given distance is $-\gamma(d)/2$ to $\gamma(d)/2$.
Using the \emph{proper} integration bounds, i.e., integrating from $d$ to $1$ for the radial limits, we leads to a large number of cases that need to be handled separately.
We, therefore, chose to follow~\cite{DBLP:journals/tog/WinchenbachA020} and integrate from $0$ to $1$, leading to complex-valued integrals for any part of the integral that would be \emph{outside} of the segment.
Note that for most terms it is necessary to treat $d=0$ as a special case of integrating angular-wise from $-\pi/2$ to $\pi/2$ and radially from $0$ to $1$, as the integral terms tend to involve divisions by $d$.
Overall this yields
\begin{equation}
\label{eqn:segment}
    \mathcal{P}^\mathcal{T}(d) = \begin{cases}
    \int_0^1\int_{\gamma(d)/2}^{\gamma(d)/2}g(r,\theta)d\theta dr, & d \neq 0\\
    \int_0^1\int_{-\pi/2}^{\pi/2} \tilde{k}(r,\theta)d\theta dr, & \text{else},
    \end{cases}
\end{equation}
with the signed distance $d$ of the origin $\mathbf{0}$ to the intersecting edge $\overline{\mathbf{\nu}_1\mathbf{\nu}_2}$ is given as
\begin{equation}
\label{eqn:dist}
    d = \frac{|(x_2-x_1)(y_1)-(x_1)(y_2-y_1)|}{\sqrt{(x_2-x_1)^2+(y_2-y_1)^2}}\cdot\begin{cases}-1,& \mathbf{0}\text{ inside}\\
    1,&\text{else}.
    \end{cases}
\end{equation}

To ensure the required spatial configuration, we first evaluate the intersections of the edge that forms the segment, i.e., the intersection of $\overline{\mathbf{\nu}_1\mathbf{\nu}_2}$ with a circle of radius $1$ centered at $\mathbf{0}$, which yields two intersection points $\mathbf{x}^p$ and $\mathbf{x}^n$.
Based on these positions, we then determine their angle of intersections, i.e., $\theta_1 = \operatorname{atan2}(x^p_y,x^p_x)$ and $\theta_2 = \operatorname{atan2}(x^n_y,x^n_x)$, and a third angle using the signed distance $\theta=\operatorname{acos}(d)$, which yields two transformation matrices $R_{\theta-\theta_1}$ and $R_{\theta-\theta_2}$. 
To simplify the integration we now ensure that the intersecting edge is orthogonal, relative to the x-axis, by checking the x-coordinate of the transformed points $\mathbf{\nu}_1$ and $\mathbf{\nu}_2$ for both transformation matrices and choosing the transformation matrix $C$ that ensures that $[C\mathbf{\nu}_{1}]_x=[C\mathbf{\nu}_{2}]_x=d$.
Accordingly, we also need to evaluate a transformed reference triangle $\mathcal{T}^\prime=C\mathcal{T}$.
Algorithm~\ref{algo:segment} describes this overall process.

\DontPrintSemicolon
\begin{algorithm}[t]
    \SetAlgoLined
    \textbf{Function} planarIntegral($\mathbf{\nu}_0$, $\mathbf{\nu}_1$, $\mathbf{\nu}_2$, $t_0$, $t_1$, $\mathcal{T}$, $f$):\;
    \Indp
    \textbf{Calculate} signed distance of $\overline{t_0t_1}$ using~\cref{eqn:dist}\;
    $\text{hit}, \lambda_1, \lambda_2 \leftarrow \operatorname{lci}(t_0, t_1-t_0,\mathbf{0},1)$\;
    $a \leftarrow t_0 + \lambda_1 (t_1 - t_0),\;b \leftarrow t_0 + \lambda_2 (t_1 - t_0)$\;
    $\theta \leftarrow\arccos d, \theta_1 = \operatorname{atan2}(a_x, a_y), \theta_2 = \operatorname{atan2}(b_x,b_y)$\;
    \textbf{If} $|(R_{\theta-\theta_1}a)_x - (R_{\theta-\theta_1}b)_x| \leq |(R_{\theta-\theta_1}a)_y - (R_{\theta-\theta_1}b)_y|$:\;
    \Indp
        $C = R_{\theta-\theta_1}$\;
    \Indm
    \textbf{Else}:
    \Indp
        $C = R_{\theta-\theta_2}$\;
    \Indm
    $\mathcal{T}^\prime\leftarrow C\mathcal{T}$\;
    \textbf{Calculate} $\tau^f_1,\tau^f_2,\tau^f_3$ for $\mathcal{T}^\prime$; see~\cref{eqn:taus}\;
    \textbf{Return} $\mathcal{P}^{\mathcal{T}^\prime}(d), C^{-1}\cdot\left[\mathcal{P}_x^{\mathcal{T}^\prime}(d),\mathcal{P}_y^{\mathcal{T}^\prime}(d)\right]^T$
    \Indm
\caption{Procedure to evaluate integrals over segments. This procedure requires as an input the edge that forms the segment described through two vertices $t_0$ and $t_1$, as well as a triangle consisting of vertices $v_0$, $v_1$ and $v_2$ where, usually, $v_1 = t_0, v_2 = t_1$. Furthermore, the reference triangle $\mathcal{T}$ and the scalar field $f$ are required. The procedure returns the integral over the region and the gradient. $\operatorname{lci}(\mathbf{o}, \mathbf{d}, \mathbf{c} ,r)$ is a function that evaluates the intersection of a line originating from $\mathbf{o}$ in direction $\mathbf{d}$ with a circle at location $\mathbf{c}$ with radius $r$ and returns an integer describing the number of intersections and the distances along $\mathbf{d}$ at which the intersections occurred.}
\label{algo:segment}
\end{algorithm}

\subsection{Integrals Over Triangle Stubs}
\label{sec:triangles:stubs}

A triangle stub is a triangle where one vertex $\mathbf{\nu}_0$, with an associated angle $\gamma$, is at the origin $\mathbf{0}$, a second vertex $\mathbf{\nu}_1$, with an associated angle $\beta$, is on the x-axis and furthest from the origin relative to the other vertices, and the last vertex $\mathbf{\nu}_2$, with an associated angle $\alpha$ is inside the support domain and above the x-axis. Initially we limit the set of triangles to those where $\alpha\leq\frac{\pi}{2}$, as this case significantly simplifies the integral. In this specific case we can split the integral over the triangle stub into first integrating radially from $0$ to $|\mathbf{\nu}_2|$ using a sector with angle $\gamma$, i.e., by utilizing $\mathcal{S}^\mathcal{T}(\gamma,|\mathbf{\nu}_2|)$, and then integrating radially from $r_\text{min}=|\mathbf{\nu}_2|$ to $r_{\text{max}}=\min\{|\mathbf{\nu}_1|,1\}$ with an angle $\omega(r)$ based on the intersection of a circle of radius $r$ with $\overline{\mathbf{\nu}_1\mathbf{\nu}_2}$. For this term we first construct a parametric line connecting $\mathbf{\nu}_1$ and $\mathbf{\nu}_2$ as $g(x) = l\tan\beta - x \tan\beta$ and a line through the origin $h(x)=x\tan\omega$, with an unknown angle $\omega$, which should intersect $g(x)$, i.e., we can solve $h(x) = g(x)$ for omega to find $\omega$ in terms of $x$ as
\begin{equation}
\label{eqn:triproblem}
    \omega(x) = \arctan\left[\frac{(l-x)\tan\beta}{x}\right]+n\pi, n\in\mathbb{Z}.
\end{equation}
However, we are interested in $\omega$ in terms of $r$. As we want the intersection to be $r$ away from the origin, we can find an orthogonal triangle with sides $x$, $y$, and $r$ with $y$ being either $g(x)$ or $h(x)$. Consequently, as $x$ is the adjacent and $r$ the hypotenuse, this yields $\cos \omega = \frac{x}{r}$, with $x = r\cos\omega$. Using~\cref{eqn:triproblem} and solving for $\omega$ then yields two solutions
\begin{equation}
    \omega_1(r) = \arcsin\left[\frac{l\sin\beta}{r}\right]-\beta +2\pi n, \omega_2(r) = -\arcsin\left[\frac{l\sin\beta}{r}\right]-\beta +2\pi n + \pi,
\end{equation}
with $n\in\mathbb{Z}$, where we utilize $\omega_1(r)$ with $n=0$, i.e., 
\begin{equation}
    \omega(r) =  \arcsin\left[\frac{l\sin\beta}{r}\right]-\beta.
\end{equation}
Accordingly, this yields an equation $\mathcal{R}$ as %the integral bounds radially are from $d = |\mathbf{\nu}_2|$ to $l = \min\{|\mathbf{\nu}_1|,1\}$ and angularily from $0$ to $\omega(r)$, and we find a term $\mathcal{R}$ as
\begin{equation}
    \mathcal{R}^\mathcal{T}(d,l,\beta,\gamma) = \mathcal{S}^T(\gamma,|\mathbf{\nu}_2|) + \int_{d}^l \int_0^{\omega(r)}g(r,\theta)d\theta dr,
\end{equation}
with analogous terms as before for $\tilde{k}_x$ and $\tilde{k}_y$. 
Given an arbitrary triangle consisting of vertices $\mathbf{\nu}_0$, at the origin, and two other vertices $\mathbf{\nu}_1$ and $\mathbf{\nu}_2$, with at least one of them within the support domain and with the correct range of angles, we first rotate the coordinate system such that the edge connecting $\mathbf{\nu}_0$ with the furthest vertex $\mathbf{\nu}_f$ is in line with the x-axis and then flipping the coordinate system if the closer vertex $\mathbf{\nu}_c$ is below the x-axis. This can be done using an improper rotation matrix $T$, which also needs to be applied to $\mathcal{T}$ to yield a transformed triangle $\mathcal{T}^T$, and yields the following function
\begin{equation}
    \mathcal{R}^\mathcal{T}(\mathbf{\nu}_0,\mathbf{\nu}_1,\mathbf{\nu}_2) = \mathcal{R}^{\mathcal{T}^T}(|\mathbf{\nu}_c|,|\mathbf{\nu}_f|,\beta,\gamma),
    \left[\begin{array}{c}
        \mathcal{R}_x^\mathcal{T}(\mathbf{\nu}_0,\mathbf{\nu}_1,\mathbf{\nu}_2) \\
        \mathcal{R}_y^\mathcal{T}(\mathbf{\nu}_0,\mathbf{\nu}_1,\mathbf{\nu}_2) 
    \end{array}\right]= T^{-1}\left[\begin{array}{c}
        \mathcal{R}_x^{\mathcal{T}^T}(|\mathbf{\nu}_c|,|\mathbf{\nu}_f|,\beta,\gamma) \\
        \mathcal{R}_y^{\mathcal{T}^T}(|\mathbf{\nu}_c|,|\mathbf{\nu}_f|,\beta,\gamma) 
    \end{array}\right], 
\end{equation}
with $|\mathbf{\nu}_c|,|\mathbf{\nu}_f| $ and $\beta$ as above. This can then be utilized to expand the range of acceptable triangles to any triangle with $\mathbf{\nu}_0 = \mathbf{0} \land  |\mathbf{\nu}_2| \leq 1 \land |\mathbf{\nu}_1| > |\mathbf{\nu}_2|$, by subdividing the triangle into two sub-triangles where the greatest angle is $\pi/2$. % see Fig.~\ref{fig:arbstub}. 
This can be done straight forwardly by determining the closest point $\mathbf{\nu}_c$ on the edge connecting $\overline{\mathbf{\nu}_1\mathbf{\nu}_2}$, using standard linear algebra, and evaluating $\mathcal{R}$ for a triangle $\mathbf{\nu}_0, \mathbf{\nu}_c, \mathbf{\nu}_1$ and a triangle $\mathbf{\nu}_0, \mathbf{\nu}_c,\mathbf{\nu}_2$ to yield
\begin{equation}
    \mathcal{G}^\mathcal{T}(\mathbf{\nu}_0,\mathbf{\nu}_1,\mathbf{\nu}_2)=\begin{cases}
        \mathcal{R}^\mathcal{T}(\mathbf{\nu}_0,\mathbf{\nu}_1,\mathbf{\nu}_2),&\gamma \leq \frac{\pi}{2},\\
        \mathcal{R}^\mathcal{T}(\mathbf{\nu}_0,\mathbf{\nu}_c,\mathbf{\nu}_1) + \mathcal{R}^\mathcal{T}(\mathbf{\nu}_0,\mathbf{\nu}_c,\mathbf{\nu}_2), &\text{ else},
    \end{cases}
\end{equation}
with analogous terms for $\mathcal{G}_x^\mathcal{T}$ and $\mathcal{G}_y^\mathcal{T}$, see Algorithm~\ref{algo:stub}.

\DontPrintSemicolon
\begin{algorithm}[t]
    \SetAlgoLined
    
    \textbf{Function} triangularIntegral($\mathbf{\nu}_1$, $\mathbf{\nu}_2$, $\mathcal{T}$, $f$):\;
    \Indp
    \textbf{If} $|\mathbf{\nu}_1| < |\mathbf{\nu}_2|$:$\mathbf{\nu}_1,\mathbf{\nu}_2\leftarrow\mathbf{\nu}_2,\mathbf{\nu}_1$\;
    $\alpha\leftarrow\arccos\frac{(\mathbf{\nu}_1-\mathbf{\nu}_2)\cdot(\mathbf{\nu}_0-\mathbf{\nu}_2)}{|\mathbf{\nu}_1-\mathbf{\nu}_2||\mathbf{\nu}_0-\mathbf{\nu}_2|}, \gamma\leftarrow\arccos\frac{(\mathbf{\nu}_1-\mathbf{\nu}_0)\cdot(\mathbf{\nu}_2-\mathbf{\nu}_0)}{|\mathbf{\nu}_1-\mathbf{\nu}_0||\mathbf{\nu}_2-\mathbf{\nu}_0|}$\;
    \textbf{If} $\gamma \leq \pi/2 \land \alpha > \pi/2$:\;
    \Indp
    \textbf{Return} integrateR($\mathbf{\nu}_1$, $\mathbf{\nu}_2$, $\mathcal{T}$, $f$)\;
    \Indm
    \textbf{Else}:\;
    \Indp
    $\mathbf{\nu}_c\leftarrow\operatorname{closestPoint}(\mathbf{0},\overline{\mathbf{\nu}_1\mathbf{\nu}_2})$\;
    $a,\nabla a = $ integrateR($\mathbf{\nu}_1$, $\mathbf{\nu}_c$, $\mathcal{T}$, $f$)\;
    $b,\nabla b = $ integrateR($\mathbf{\nu}_2$, $\mathbf{\nu}_c$, $\mathcal{T}$, $f$)\;
    \textbf{Return} $a + b, \nabla a + \nabla b$\;
    \Indm
    \Indm
    
    \textbf{Function} integrateR($\mathbf{\nu}_1$, $\mathbf{\nu}_2$, $\mathcal{T}$, $f$):\;
    \Indp
    $\mathbf{\nu}_0 \leftarrow \mathbf{0}$\;
    $l \leftarrow |\mathbf{\nu}_1|, d\leftarrow |\mathbf{\nu}_2|, \theta - \operatorname{atan2}(\mathbf{\nu}_1^x,\mathbf{\nu}_1^y)$\;
    $s\leftarrow \operatorname{sgn}(\mathbf{\nu}_2^x\sin(-\theta)+\mathbf{\nu}_2^t\cos(-\theta)) >= 0 ? 1 : -1$\;
    $\beta\leftarrow\arccos\frac{(\mathbf{\nu}_2-\mathbf{\nu}_1)\cdot(\mathbf{\nu}_0-\mathbf{\nu}_1)}{|\mathbf{\nu}_2-\mathbf{\nu}_1||\mathbf{\nu}_0-\mathbf{\nu}_1|}, \gamma\leftarrow\arccos\frac{(\mathbf{\nu}_1-\mathbf{\nu}_0)\cdot(\mathbf{\nu}_2-\mathbf{\nu}_0)}{|\mathbf{\nu}_1-\mathbf{\nu}_0||\mathbf{\nu}_2-\mathbf{\nu}_0|}$\;
    $T\leftarrow\left[\begin{array}{ccc}
         \cos\theta&\sin\theta&0  \\
         -s\sin\theta&s\cos\theta&0\\
         0&0&1
    \end{array}\right]$\;
    \textbf{Return} $\mathcal{R}^\mathcal{T}(l,d,\beta,\gamma)$\;
    \Indm

\caption{Procedure to evaluate integrals over triangle stubs. This procedure requires two vertices $\mathbf{\nu}_1$ and $\mathbf{\nu}_2$ of a triangle as the third vertex $\mathbf{\nu}_0$ is assumed to lie at the origin. Furthermore, the reference triangle $\mathcal{T}$ and the scalar field $f$ are required. The procedure returns the integral over the region and the gradient.}
\label{algo:stub}
\end{algorithm}

\subsection{Integrals Over Arbitrary Wedges}
\label{sec:triangles:wedges}

An arbitrary wedge is any area that results from the intersection of a triangle with a circle where $0$ or $1$ vertices are inside of the circle.
Note that this also includes cases where there is no overlap, i.e., the integral is $0$, or the triangle completely overlaps the circle, i.e., the integral is based on the integral over the entire support domain; however, these cases can readily be detected and handled.
Furthermore, there being only one edge intersecting the circle, with no vertices inside, is identical to the segment case discussed prior.
By convention, we assume that either no vertex is inside of the circle or that $\mathbf{\nu}_0$ is in the circle, with the triangle being integrated having positive area, which can both be ensured by reordering the vertex indices accordingly.
To handle the special cases, we first evaluate the intersections of all edges with the circle and if no edge intersects the circle we either return $0$, if the origin is not in the triangle, an integral over the entire support domain if the origin is in the triangle and an integral over a segment if only a single edge intersects the support domain.
Consequently, we only need to handle cases where at least $2$ edges intersect the support domain and at most $1$ vertex, $\mathbf{\nu}_0$, is inside of the support domain.

Accordingly, we can find the intersections of the adjacent edges, of $\mathbf{\nu}_0$, with the circle as $\lambda_n^1,\lambda_p^1,\lambda_n^2$ and $\lambda_p^2$ with $\lambda_n < \lambda^p$.
As $\mathbf{\nu}_1$ and $\mathbf{\nu}_2$ both lie outside the circle, $\lambda_p$ will always be on the corresponding edge and lie on the circle, whereas $\lambda_n$ will only be positive if $\mathbf{\nu}_0$ was also outside the circle. 
Next, we determine the two intersecting points $\mathbf{\sigma}^1 = \mathbf{\nu}_0 + \lambda_p^1 \left[\mathbf{\nu}_1 - \mathbf{\nu}_0\right]$ and $\mathbf{\sigma}^2 = \mathbf{\nu}_0 + \lambda_p^2 \left[\mathbf{\nu}_2 - \mathbf{\nu}_0\right]$, which both lie on the circle, with $|\mathbf{\sigma}^1|=|\mathbf{\sigma}_2|=1$.
Using these intersecting points, we can then split the integral over the original triangle $\mathcal{T}$ as an integral over the arc between $\mathbf{\sigma}^1$ and $\mathbf{\sigma}_2$ and the integral over the generic triangle stubs $\mathcal{T}_1 = \{\mathbf{0},\mathbf{\nu}_r,\mathbf{\sigma}^1\}$ and $\mathcal{T}_2 = \{\mathbf{0},\mathbf{\sigma}^2,\mathbf{\nu}_r\}$, i.e., we find
\begin{equation}
    \iint_\mathcal{T}\tilde{k}(r,\theta)d\theta dr =
    C^\mathcal{T}(\angle\sigma^1\sigma^2,1) + G^\mathcal{T}(\mathbf{0},\mathbf{\nu}_0,\sigma^1)+ G^\mathcal{T}(\mathbf{0},\sigma^2,\mathbf{\nu}_0).
\end{equation}

However, this argumentation only works if the origin $\mathbf{0}$ is within the triangle. To handle cases where the origin is not within the triangle $\mathcal{T}$, we multiply each integral result by the sign of the area that is being used for the integration, i.e., we instead use 
\begin{equation}
\label{eqn:intermediateWedge}
\begin{split}
    \iint_\mathcal{T}\tilde{k}(r,\theta)d\theta dr =&
    \operatorname{sgn}\left(A(\{\mathbf{0},\sigma^1,\sigma^2\})\right)C^\mathcal{T}(\angle\sigma^1\sigma^2,1) \\ +&\operatorname{sgn}\left(A\{\mathbf{0},\mathbf{\nu}_0,\sigma^1\})\right)G^{\mathcal{T}}(\mathbf{0},\mathbf{\nu}_0,\sigma^1)
    +\operatorname{sgn}\left(A\{\mathbf{0},\sigma^2,\mathbf{\nu}_0\})\right)G^{\mathcal{T}}(\mathbf{0},\sigma^2,\mathbf{\nu}_0).
    \end{split}
\end{equation}

If the edge opposite of the reference vertex $\mathbf{\nu}_r$, i.e., $\overline{\mathbf{\nu}_1\mathbf{\nu}_2}$ intersects the circle then we find two intersection points along this line $\lambda_n^3$ and $\lambda_p^3$, which describe a segment that was included, when it shouldn't have, during the integration over the cone from $\sigma^1$ to $\sigma^2$ and, accordingly, this region needs to be removed from the integral. 
This region can be described as a segment that has been solved before, where we only need to consider how we can formulate the problem at hand in a way that makes reusing the prior results possible.
To achieve this, we first determine a new vertex $\mathbf{\nu}_c$ that lies outside of the circle and on the opposite side, relative to the origin, of $\overline{\mathbf{\nu}_1\mathbf{\nu}_2}$, and using $\mathcal{P}^\mathcal{T}$ with the triangle $\mathbf{\nu}_c, \mathbf{\nu}_1, \mathbf{\nu}_2$.
We then determine the sign of the signed area of the triangle $\mathbf{\nu}_c, \mathbf{\nu}_1, \mathbf{\nu}_2$ and subtract $\operatorname{sgn}(A(\{\mathbf{\nu}_c, \mathbf{\nu}_1, \mathbf{\nu}_2\}))\mathcal{P}^\mathcal{T}(\mathbf{\nu}_c, \mathbf{\nu}_1, \mathbf{\nu}_2)$ from the integral over the cone, i.e., $C^\mathcal{T}(\angle\sigma^1\sigma^2,1)$ in~\cref{eqn:intermediateWedge}.
If the signed area of the sector described by $\sigma_1$ and $\sigma_2$ was positive then the integral over the planar region has to be flipped, i.e., the integral has to be evaluated over everything but this region, which can readily be achieved by calculating $\mathcal{S}^\mathcal{T}(1) - \mathcal{P}^\mathcal{T}(\mathbf{\nu}_c, \mathbf{\nu}_1, \mathbf{\nu}_2)$.
Accordingly, we can find a function $\mathcal{W}^\mathcal{T}(\mathbf{\nu}_0,\mathbf{\nu}_1,\mathbf{\nu}_2)$ that works for all potential cases with $0$ or $1$ vertex inside the circle and $0$ to $3$ edges intersecting the circle; see Algorithm~\ref{algo:Wedge}.

\DontPrintSemicolon
\begin{algorithm}[t]
\SetAlgoLined
\textbf{Function} integrateWedge($\mathbf{\nu}_0$, $\mathbf{\nu}_1$, $\mathbf{\nu}_2$, $\mathcal{T}$, $f$):\;
\Indp
    \textbf{If} no edges intersect:\;
    \Indp
    \textbf{If} $\mathbf{O} \in \text{Triangle}\{\mathbf{\nu}_0,\mathbf{\nu}_1,\mathbf{\nu}_2\}$: \textbf{Return} $0, [0,0]^T$\;
    \textbf{Else}: \textbf{Return} $\mathcal{S}(1)$\;
    \Indm
    \textbf{If} one edge intersects: \textbf{Return} segment integral\;
    \;
    $h_1, \lambda_n^1, \lambda_p^1 \leftarrow \operatorname{lci}(\mathbf{\nu}_0,\mathbf{\nu}_1-\mathbf{\nu}_0,\mathbf{0},1)$\;$ \lambda_1 \leftarrow \lambda_p^1 \leq 1 ? \lambda_p^1 : \lambda_n^1$\;
    $h_2, \lambda_n^2, \lambda_p^2 \leftarrow \operatorname{lci}(\mathbf{\nu}_0,\mathbf{\nu}_2-\mathbf{\nu}_0,\mathbf{0},1)$\;$\lambda_2 \leftarrow \lambda_p^2 \leq 1 ? \lambda_p^2 : \lambda_n^2$\;
    
    $\sigma_1 \leftarrow \mathbf{\nu}_0 + \lambda_1 (\mathbf{\nu}_1-\mathbf{\nu}_0), \sigma_2 \leftarrow \mathbf{\nu}_0 + \lambda_2 (\mathbf{\nu}_2-\mathbf{\nu}_0)$\;
    $q \leftarrow \operatorname{Area}(\text{tri})\geq 0 ? 1 : 0, p \leftarrow \mathbf{0}\in\text{tri}? \text{True}:\text{False}$\;
    \;
    
    $a , \nabla a \leftarrow \operatorname{triangularIntegral}(\mathbf{\nu}_0, \sigma_1,\mathcal{T},f)$\; $s_a = \operatorname{sgn}(\operatorname{Area}(\{\mathbf{0},\mathbf{\nu}_0,\sigma_1\}))$\;
    $b , \nabla b \leftarrow \operatorname{triangularIntegral}(\sigma_2, \mathbf{\nu}_0,\mathcal{T},f)$\;$s_b = \operatorname{sgn}(\operatorname{Area}(\{\mathbf{0},\sigma_2,\mathbf{\nu}_0\}))$\;
    $c,\nabla c \leftarrow \operatorname{sectorIntegral}(\sigma_1,\sigma_2,1,\mathcal{T},f)$\;
    
    \textbf{If} third edge intersects:\;%$h=\text{True} \land (l_1 \geq 0 \land l_1 \leq 1 \land l_2 \geq 0 \land l_2 \leq 1)$:\;
    \Indp
    $s,\nabla s \leftarrow \operatorname{sphericalIntegral}(1,\mathcal{T},f)$\;
    $\mathbf{\nu}_c \leftarrow \left[\mathbf{\nu}_1^x - 2 \left(\mathbf{\nu}_2^y-\mathbf{\nu}_1^y\right),\mathbf{\nu}_1^y + 2 \left(\mathbf{\nu}_2^x-\mathbf{\nu}_1^x\right)\right]^T$\;
    \textbf{If} $\mathbf{\nu}_c$ on same side as $\mathbf{0}$ relative to $\overline{\mathbf{\nu}_1\mathbf{\nu}_2}$:\;
    \Indp
        $\mathbf{\nu}_c \leftarrow \left[\mathbf{\nu}_1^x + 2 \left(\mathbf{\nu}_2^y-\mathbf{\nu}_1^y\right),\mathbf{\nu}_1^y - 2 \left(\mathbf{\nu}_2^x-\mathbf{\nu}_1^x\right)\right]^T$\;
    \Indm
    $p,\nabla p \leftarrow \operatorname{planarIntegral}(\mathbf{\nu}_c,\mathbf{\nu}_1,\mathbf{\nu}_2,\mathbf{\nu}_2,\mathbf{\nu}_1)$\;
    \textbf{If} $q \geq 0$ : $p,\nabla p \leftarrow s - p, \nabla s - \nabla p$\; 
    $c,\nabla c \leftarrow c - s_t p, \nabla c - s_t \nabla p$\;
    \Indm
    $p \leftarrow q c + s_a a + s_b b$\;
    $\nabla p \leftarrow q \nabla c + s_a \nabla a + s_b \nabla b$\;
    
    \textbf{Return} $p,\nabla p $\;

\Indm
    
\caption{Procedure to evaluate integrals over arbitrary wedges. This procedure requires three vertices $\mathbf{\nu}_0, \mathbf{\nu}_1$ and $\mathbf{\nu}_2$ of a triangle, the reference triangle $\mathcal{T}$ and the scalar field $f$. The procedure returns the integral over the region and the gradient, with $\operatorname{lci}$ defined as in~\ref{algo:sector}}
\label{algo:Wedge}
\end{algorithm}

\subsection{Integrals Over Triangles}
\label{sec:triangles:triangles}

Considering the previous section, the only remaining cases that need to be handled are cases in which either $2$ or $3$ vertices are within the circle. These cases can be constructed straight-forwardly using the arbitrary wedge function $\mathcal{W}^\mathcal{T}(\mathbf{\nu}_0,\mathbf{\nu}_1,\mathbf{\nu}_2)$ by constructing arbitrary wedges that resemble these cases. In the case of two vertices inside the circle, i.e., $\mathbf{\nu}_0$ and $\mathbf{\nu}_1$ we can construct a point $\mathbf{\nu}_x$ that lies along the extended edge $\overline{\mathbf{\nu}_0\mathbf{\nu}_1}$ and is outside of the circle, e.g., $\mathbf{\nu}_x = \mathbf{\nu}_0 + 3 \frac{\mathbf{\nu}_1-\mathbf{\nu}_0}{|\mathbf{\nu}_1-\mathbf{\nu}_0|}$, which yields two triangles $\mathcal{T}_1 = \{\mathbf{\nu}_0,\mathbf{\nu}_x,\mathbf{\nu}_1\}$ and $\mathcal{T}_2=\{\mathbf{\nu}_0,\mathbf{\nu}_x,\mathbf{\nu}_2\}$, which are wedge shaped, where the overall integral then is \begin{equation}
    T_2^\mathcal{T}(\mathbf{\nu}_0,\mathbf{\nu}_1,\mathbf{\nu}_2) = \mathcal{W}^\mathcal{T}(\mathbf{\nu}_0,\mathbf{\nu}_x,\mathbf{\nu}_1) - \mathcal{W}^\mathcal{T}(\mathbf{\nu}_0,\mathbf{\nu}_x,\mathbf{\nu}_2).
\end{equation}
Note that analogous cases can be constructed if $\mathbf{\nu}_1$ and $\mathbf{\nu}_2$ are inside the circle and if $\mathbf{\nu}_0$ and $\mathbf{\nu}_2$ are inside the circle.
In case of three vertices inside of the circle we construct a similar point $\mathbf{\nu}_x^1$ along the extended edge $\overline{\mathbf{\nu}_0\mathbf{\nu}_1}$ and $\mathbf{\nu}_x^2$ along the extended edge $\overline{\mathbf{\nu}_0\mathbf{\nu}_2}$. The integral can then be evaluated as 
\begin{equation}
     T_3^\mathcal{T}(\mathbf{\nu}_0,\mathbf{\nu}_1,\mathbf{\nu}_2) = \mathcal{W}^\mathcal{T}(\mathbf{\nu}_0,\mathbf{\nu}_x^1,\mathbf{\nu}_x^2) - T_2^\mathcal{T}(\mathbf{\nu}_1,\mathbf{\nu}_2,\mathbf{\nu}_x^1) -\mathcal{W}^\mathcal{T}(\mathbf{\nu}_2,\nu^1_x,\nu^2_x)
\end{equation}

\subsection{Overall integral}
Based on the results for wedges $\mathcal{W}$, i.e., triangles with $0$ or $1$ vertex inside of the triangle, and the results for $2$, $T_2$, and $3$, $T_3$, vertices in the circle the overall integral for an arbitrary triangle $\mathcal{T}=\{\mathbf{\nu}_0,\mathbf{\nu}_1,\mathbf{\nu}_2\}$can be evaluated based on the number of vertices inside the triangle $\#\nu$ as
\begin{equation}
    \iint_\mathcal{T}g(r,\theta)d\theta dr = \begin{cases}
    \mathcal{W}^\mathcal{T}(\mathbf{\nu}_0,\mathbf{\nu}_1,\mathbf{\nu}_2),& \#\nu \leq 1,\\
    T_2^\mathcal{T}(\mathbf{\nu}_0,\mathbf{\nu}_1,\mathbf{\nu}_2),& \#\nu = 2,\\
    T_3^\mathcal{T}(\mathbf{\nu}_0,\mathbf{\nu}_1,\mathbf{\nu}_2),& \#\nu = 3.\\
    \end{cases}
\end{equation}
Finally, it is important to note that utilizing gradient terms that are not using the straightforward formulation, see~\cref{eqn:basicdiff}, the methodology is unchanged but requires calculating the weights $\tau_i^f$ using modified terms that were discussed prior at the end of Section~\ref{sec:integrals}.

\section{Closed Form Solutions}
\label{sec:integrals}

Based on the previous observations in Section~\ref{sec:barycentric}, we can reduce the set of functions that we need to integrate to the following set of three terms
\begin{align}
    g_n^{0,0} &= x^n,\\
    g_n^{a,0} &= x^n\cos\left(a\theta\right),\\
    g_n^{a,b} &= x^n\cos\left(a\theta\right)\sin\left(b\theta\right).\\
\end{align}

We can further simplify this result using a trigonometric transformation for the third term to remove the cosine term, i.e., we can use $cos(a)sin(b) = 1/2 sin(a+b) - 1/2sin(a-b)$, which in turn yields
\begin{equation}
    \frac{1}{2}\sin\left((a+b)\theta\right)x^n - \frac{1}{2}\sin\left((a-b)\theta\right)x^n.
\end{equation}

However there are two important notes about this function when integrating it as (a) if $a-b=0$ the second term should not be evaluated, i.e., even though we can compute an integral that is valid for this term, it would not lead to meaningful results and should be avoided, and (b) both terms are more general forms of $\sin(b\theta)$.
This means we only need to integrate
\begin{align}
    f_p(x) &= x^n\\
    f_c(x,\theta) &= x^n\cos\left(a\theta\right)\\
    f_s(x,\theta) &= x^n\sin\left(b\theta\right).
\end{align} 

Before discussing the actual solutions for this integral, we will first discuss more general closed form terms that will be at the core of the integral terms in Section~\ref{sec:integrals:fundamental}.
After that we will discuss integrals for conical segments briefly, see Sec.~\ref{sec:integrals:cones}, followed by a more in-depth derivation of the terms for segments, in Sec.~\ref{sec:integrals:segments}, and stubs in Sec.~\ref{sec:integrals:stubs}.

\subsection{Fundamental Closed Form Terms}
\label{sec:integrals:fundamental}

The first general integral we will consider is for the term 
\begin{equation}
\int x^n\sqrt{1-\frac{d^2}{x^2}}^m dx,
\end{equation}
which occurs in the majority of terms for segments and stubs.
This term has a special solution for most inputs we consider and reads
\begin{equation}
\int x^n\sqrt{1-\frac{d^2}{x^2}}^m dx = -\frac{x^{n-m+1}d^{m-1}\sqrt{x^2-d^2}}{\sqrt{1-\frac{x^2}{d^2}}(m-1-n)}~_2F_1\left(-\frac{m}{2},\frac{n - 2(m-1)}{2},\frac{2+n-2(m-1)}{2};\frac{x^2}{d^2}\right). 
\end{equation}
There are several special cases where this function becomes ill-defined, e.g., for $n=1\land m=2$.
However, the most crucial challenging term is $\frac{\sqrt{x^2-d^2}}{\sqrt{1-\frac{x^2}{d^2}}}$ for $x=d$, which is a common lower bound, becomes $0/0$.
However, by looking at the limits we find that
\begin{equation}
\lim_{x->d^-} \frac{\sqrt{x^2-d^2}}{\sqrt{1-\frac{x^2}{d^2}}} = \frac{\sqrt{-d}}{\sqrt{\frac{1}{d}}}=\sqrt{-d}\sqrt{d} =i d;\;\lim_{x->d^+} \frac{\sqrt{x^2-d^2}}{\sqrt{1-\frac{x^2}{d^2}}} = \frac{\sqrt{d}}{\sqrt{-\frac{1}{d}}}=-\sqrt{-\frac{1}{d}} d^{3/2}=-id.
\end{equation}
As this term specifically occurs on the lower end of the integral bounds, we pick the right limit of the fraction by evaluating $x+\epsilon$, with an arbitrary small epsilon, if $x=d$.
There are six other special cases that we need to treat separately, and are given as
\begin{align}
n=1,m=0&\rightarrow log(x),\\
n=1,m\neq0&\rightarrow \frac{1}{n+1}x^{n+1},\\
n=0,m=1&\rightarrow \frac{x\sqrt{-d^2+x^2}-\frac{\operatorname{atan}\left(\sqrt{-d^2+x^2}\right)}{x}}{|x|},\\
n=1,m=2&\rightarrow \frac{1}{2} \left(-\frac{d^4}{x^2}+x^2-4 d^2\log{x}\right)\\
n=3,m=2&\rightarrow -d^2 x^2 + \frac{x^4}{4}+d^4\log{x}\\
n\geq3,m=2&\rightarrow x^{n-3}\left(\frac{d^4}{n-3}-\frac{2d^2x^2}{n-1} + \frac{x^4}{n+1}\right)
\end{align}

Evaluating the general term, however, is not as straightforward as these special cases due to involving the $~_2F_1$, which in general is itself defined as the solution of an integral, i.e., the term is not algebraic for all values.
However, for the specific subset of terms we require (which is one of the primary limitations to integer power polynomials), this term can be simplified to a series of algebraic terms that can be evaluated.
At the core of these simplifications are a few special cases given as
\begin{equation}
    ~_2F_1(0,b,c,z) = 1,\; ~_2F_1(-1,b,c,z) = \frac{c-bz}{c},\;~_2F_1(a,b,b,z) = (1-z)^{-a},\;F(a,b,c,0)=1,
\end{equation}
which are straightforward to implement. 
Furthermore, using Gauss' contiguous relations allows one to relate any set of three $~_2F_1$ terms through simple relations if they differ only by integer multiples.
Specifically, we utilize
\begin{equation}
(c-1)(~_2F_1(a,b,c-1,z) - ~_2F_1(a,b,c,z)) = b(~_2F_1(a,b+1,c,z) - ~_2F_1(a,b,c,z)) = a(~_2F_1(a+1,b,c,z) - ~_2F_1(a,b,c,z)),
\end{equation}
which can be used to evaluate the series for arbitrary integer offsets of $a$ from $1$ and $1/2$.
Using the identities, the only remaining term that we will need is $~_2F_1(-1/2+\frac{m}{2},{1+b}{2},{b+3}{2},z)$ as all other required terms can be derived from this term through Gauss' relations.
We can then find a closed-form expression for this term for even and odd $m$, i.e., terms for which the second and third arguments are either integers or half steps away from integers.
To define these series expressions, we first define a set of polynomial factors resulting from the relation of the $~_2F_1$ term to Chebyshev polynomials as
\begin{equation}
S(n) = \binom{2n}{n}\frac{1}{2\cdot4^n},\; T(n) = \frac{1}{(2n+1)S(n)},
\end{equation}
which can then be utilized to define a central series term of the function as
\begin{equation}
F(n, c) = 2 \frac{\sin^{-1}\sqrt{c}}{\sqrt{1-c}} - \sum_{i=0}^{n/2-1} T(i) c^{i+1},
\end{equation}
where an important consideration is the definition of the complex inverse $\sin$ function.
To find consistent results, we utilize the definition used in Mathematica where $\sin^{-1}(c) = -\sin^{-1}(-c)$, even for complex-valued arguments, which stems from the definition using the integral extension, instead of using the logarithm-based extension used in Numpy.
Another important note is that this term is very sensitive to numerical errors, as the terms for larger polynomial degrees can become very large, leading to issues of elimination.
These, for the cases we consider, can be resolved using Horner's scheme to evaluate the polynomial series on the right side, but in general would need to be handled explicitly, e.g., by using an explicit series representation of the left-hand term.
Using this term, the odd Hypergeometric term then becomes
\begin{equation}
B_\text{odd}(n,c) = \sqrt{1-c}\left(\frac{n}{n+1} + S\left(\frac{n+1}{2}\right) \frac{1}{c^{n/2}} F(n,c)\right).
\end{equation}
For even order terms, we first find an equivalent series term
\begin{equation}
G(n,c) = \sum_{i=0}^{n/2+1}S(j)\sqrt{1-c} - \frac{1}{2}.
\end{equation}
This in turn leads to the even Hypergeometric terms becoming
\begin{equation}
B_\text{even}(n,c) = \sqrt{1-c} - \beta\left(\frac{n}{2}+1,\frac{n}{2}\right)\left(G(n,c) c^{-n/2}\right),
\end{equation}
with $\beta$ being the complete $\beta$ function defined using the complete $\Gamma$ function as $\beta(a,b)=\Gamma(a) \Gamma(b) / \Gamma(a + b)$.
Another integral that we require in several places is an integral of the inverse secant function, which can be found in standard integral lists and is given as
\begin{equation}
\begin{split}
\int_a^b\operatorname{sec}^{-1}\left(ax\right)x^m = -&\frac{b^{m+1}}{(m+1)^2}\left((1+m)\operatorname{sec}^{-1}\left(ab\right)+a\sqrt{1-\frac{1}{a^2b^2}}b ~_2F_1\left(1,\frac{2+m}{2},\frac{3+m}{2},a^2b^2\right)\right)\\
&+\frac{c^{m+1}}{(m+1)^2}\left((1+m)\operatorname{sec}^{-1}\left(ac\right)+a\sqrt{1-\frac{1}{a^2c^2}}c ~_2F_1\left(1,\frac{2+m}{2},\frac{3+m}{2},a^2c^2\right)\right),
\end{split}
\end{equation}
which can be evaluated as above.
The final set of integrals we require is based on the Chebyshev Polynomials of fthe irst and second kind, which canonically are given as
\begin{align}
    U_n(x) &= \sum_{k=0}^{\lfloor n/2\rfloor}(-1)^k\binom{n-k}{k}2^{n-2k} x^{n-2k},\\
    T_n(x) &= \sum_{m=0}^{\lfloor n/2\rfloor}(-1)^m\left(\binom{n-m}{m}+\binom{n-m-1}{n-2m}\right)\cdot2^{n-2m-1}\cdot x^{n-2m}.
\end{align}
When integrating these terms with respect to $x$, we can directly utilize the series expansion.
While this does not generally lend itself to finding compact expressions, it does lead to terms that are straightforward to integrate by changing the order of summation and integration.
This does potentially lead to a large number of evaluations of integrals, e.g., the square root integral described before, but still yields accurate results.
\subsection{Integrals over Cones}
\label{sec:integrals:cones}

Integrating over conical segments is the most fundamental case we need to consider and consists of solving the integral with bounds 
\begin{equation}
    \int_0^l\int_\alpha^\beta f(x,\theta)d\theta dx,
\end{equation}
Which for our three terms of concern yields
\begin{align}
\int_0^l\int_\alpha^\beta f_p(x) d\theta dx = \frac{\left( \beta\operatorname{-}\alpha\right)  {{l}^{n\operatorname{+}1}}}{n\operatorname{+}1},\\
\int_0^l\int_\alpha^\beta f_c(x,\theta) d\theta dx = \left(\frac{\sin(a\beta)}{a}-\frac{\sin(a\alpha)}{a}\right)\frac{l^{n+1}}{n+1},\\
\int_0^l\int_\alpha^\beta f_c(x,\theta) d\theta dx = \left(\frac{\cos(b\beta)}{b}-\frac{\cos(b\alpha)}{b}\right)\frac{l^{n+1}}{n+1}.
\end{align}
Note that it can be useful for efficiency to precompute certain special cases, e.g., integrals over the entire support domain, or half the support domain, in advance as these terms occur very regularly.
Otherwise, this is all that is necessary for these terms.
\subsection{Integrals over Segments}
\label{sec:integrals:segments}

For circular segments we require an integral from $-\cos^{-1}\left(\frac{d}{x}\right)$ to $\cos^{-1}\left(\frac{d}{x}\right)$, where contrary to the previous section, we will first consider $f_s$ as this integral evaluates to $0$. 
Next, we consider the integral for the cosine terms, i.e., $f_c$, which, after integrating the angular component, yields an integral of the form
\begin{equation}
\int_0^1 \frac{2}{a} sin\left(a \cos^{-1}\left(\frac{d}{x}\right)\right)x^n dx.
\end{equation}
To simplify this term, we utilize a somewhat lesser-known identity of the Chebyshev polynomials that is given as:
\begin{equation}
\sin\left(n \cos^{-1}(x)\right) = \sqrt{1 - x^2} U_{n-1}(x),
\end{equation}
which means that we can insert this term into the integral first
\begin{equation}
\frac{2}{a} \int_0^1 \sqrt{1-\frac{d^2}{x^2}}U_{a-1}\left(\frac{d}{x}\right)x^n dx,
\end{equation}
which is straightforward to integrate by inserting the series definition of the Chebyshev polynomials of the second kind and flipping the order of integration and summation, i.e., we can expand the term to yield
\begin{equation}
\frac{2}{a} \sum_{k=0}^{\lfloor \frac{a-1}{2}\rfloor}(-1)^k\binom{a-1-k}{k}2^{a-1-2k} d^{a-1-2k} \int_0^1 \sqrt{1-\frac{d^2}{x^2}}x^{2k+n+1-a} dx,
\end{equation}
for which we can use the previously described solution of the square root integral terms.
This leaves the integral over the polynomial terms $f_p$, which gives an integral of
\begin{equation}
2 x^n \cos ^{-1}\left(\frac{d}{x}\right).
\end{equation}
Integrating this term is not straight forward directly, however, we can utilize a trigonometric identity as $\cos ^{-1}\left(x\right) = \sec^{-1}\left(\frac{1}{x}\right)$, to replace the cosine term with a secant term.
For this secant term, there exists a well-known closed-form integral, which we described in Sec.~\ref{sec:integrals:fundamental}, which can be directly utilized to evaluate this integral.
Note that for all of the cases described here, we still utilize the cone segment solutions for $d=-1$, $d=0$, and $d=1$, as these terms tend to require special treatment.

\subsection{Integrals over Stubs}
\label{sec:integrals:stubs}
Moving to the stub triangles we need an integral of the form
$$
\int_l^u\int_0^{-\beta+\sin^{-1}\frac{c\sin\beta}{x}}f_n^{a,b}(x,\theta)d\theta dx
$$

For readability we can drop the $\sin\beta$ and assume the bounds to be $\sin^{-1}\left(\frac{d}{x}\right)-\beta$. %Inserting for the first term we find
Inserting $f_p$ yields the following term:
\begin{equation}
\left(\operatorname{asin}\left( \frac{\sin{d}}{x}\right) \operatorname{-}\beta\right)  {{x}^{n}} = {x}^{n}\operatorname{asin}\left( \frac{\sin{d}}{x}\right) \operatorname{-}\beta{x}^{n},
\end{equation}
where the first term can be integrated as before, i.e., using the arc secant integral term, and the second term is trivial to integrate.
Integrating $x^n \cos a\theta$ gives
\begin{equation}
  \frac{\sin{\left( a \operatorname{asin}\left( \frac{\sin{(\beta)} l}{x}\right) \operatorname{-}a \beta\right) } {{x}^{n}}}{a},  
\end{equation}
which can be simplified using the trigonometric identity $sin(a+b) = cos(a)sin(b)+sin(a)cos(b)$ to find
$$\frac{\cos(-a\beta)}{a} x^n \sin\left(a \sin^{-1}\left[\frac{d}{x}\right]\right) + \frac{\sin(-a\beta)}{a}x^n \cos\left(a \sin^{-1}\left[\frac{d}{x}\right]\right)$$

Integrating the $x^n\sin b\theta $term analogously yields
\begin{equation}
    \frac{\left( \frac{1}{b\operatorname{+}a}\operatorname{-}\frac{\cos{\left( b  \operatorname{asin}\left( \frac{\sin{\beta} l}{x}\right) \operatorname{+}-b  \beta\right) }}{b\operatorname{+}a}\right)  {{x}^{n}}}{2}
\end{equation}
Simplifying and expanding the trigonometric terms gives 
\[\frac{\left( \sin{\left( -b  \beta\right) } \sin{\left( b  \operatorname{asin}\left( \frac{\sin{(beta)} l}{x}\right) \right) }\operatorname{-}\cos{\left( -b  \beta\right) } \cos{\left( b  \operatorname{asin}\left( \frac{\sin{(\beta)} l}{x}\right) \right) }\operatorname{+}1\right)  {{x}^{n}}}{2 b }\]

For both terms, the central terms are variants of 
\begin{align}
 \cos\left(n\operatorname{asin}\left( \frac{d l}{x}\right)\right),\;\sin\left(n\operatorname{asin}\left( \frac{d l}{x}\right)\right),
\end{align}
which can be utilized using the definition of the Chebyshev series, similar to before, as 

\begin{align}
\sin a \left(\sin^{-1}\left(\frac{d}{x}\right)\right) &= x U_{a-1}\left(\sqrt{1 - \frac{d^2}{x^2}}\right)\\
\cos \left(a \sin^{-1}\left(\frac{d}{x}\right)\right) &= T_a\left(\sqrt{1 - \frac{d^2}{x^2}}\right)
\end{align}
Both of which are straightforward to solve using the square root integral term derived before, and by inserting the definition of the Chebyshev series as before, and then relying on our solution for the square root term, as all other terms are directly derived from this one.

\section{Evaluation}
\label{sec:eval}
Our proposed method offers an exact analytical and algebraic solution for evaluating boundary integrals over arbitrary 2D triangles, accommodating arbitrary-order quantities and their derivatives. 
Consequently, evaluating the accuracy of our solution is difficult outside of limited cases with known solutions.
Instead, we evaluate its accuracy and efficacy by demonstrating the convergence of established numerical quadrature schemes toward our analytical solution in four distinct evaluation setups:
First, we qualitatively compare our solution against conventional numerical quadrature rules on a single triangle. This highlights scenarios where existing methods typically introduce errors, underscoring the precision of our analytical approach.
Second, we quantitatively assess convergence using a well-defined test case: random triangulations of $[-1,1]^2$ with a known global integral and gradient for a piecewise linear quantity. This allows us to observe numerical scheme convergence towards the known solution and allows us to validate our method against a globally exact result, offering insights into its implementation and correctness.
Finally, we showcase the practical integration of our boundary integral method within a simple incompressible SPH scheme. 
While a comprehensive application study is beyond this paper's scope, this demonstration illustrates the potential for enhancing CFD simulations by incorporating our analytical boundary integral solutions.

\subsection{Qualitative Evaluation of Quadrature Schemes}
\begin{figure}
    \centering
    \includegraphics[width=0.75\linewidth]{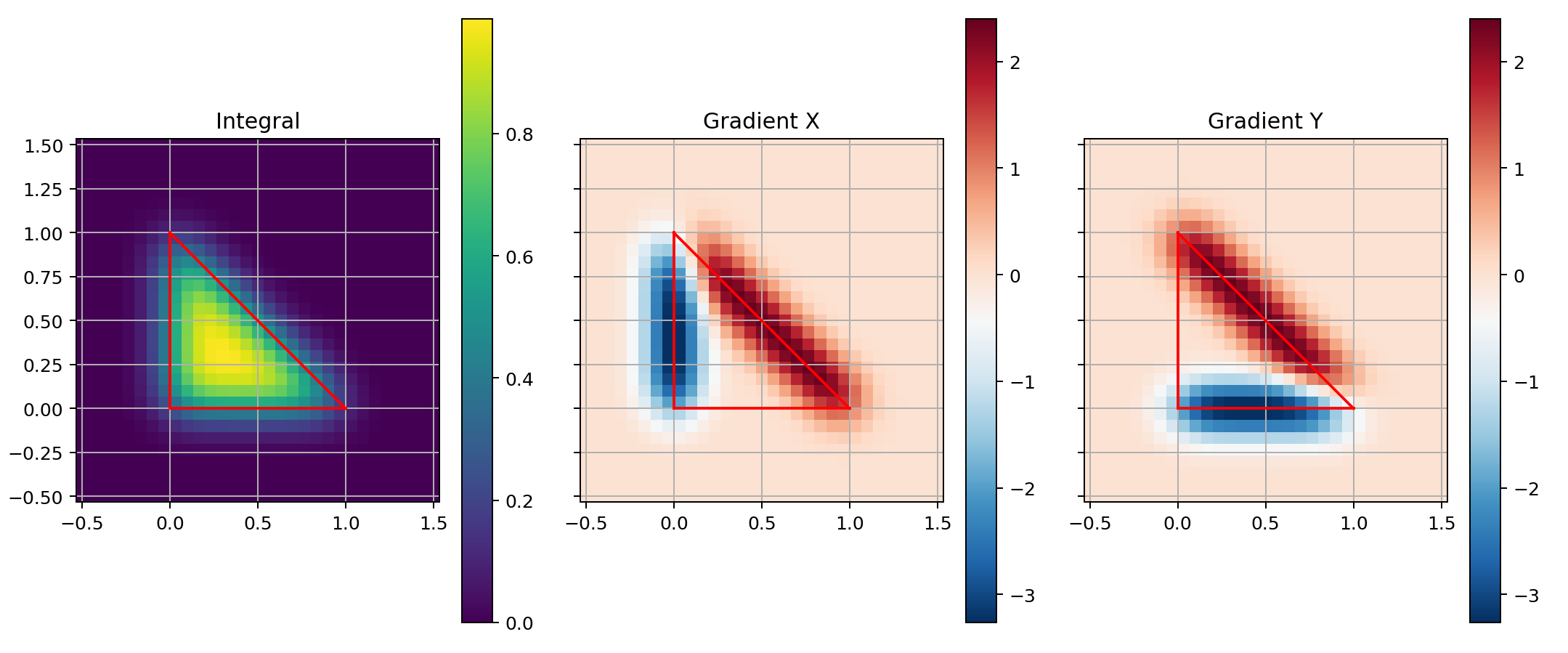}
    \caption{Exemplary evaluation of the boundary integral results for a Cartesian grid of coordinates surrounding a triangle with support $h=1$, showing the result for the integral and the gradient for a piecewise-constant boundary quantity of $1$.}
    \label{fig:smallTriangleExample}
\end{figure}

\begin{figure}
    \centering
    \includegraphics[width=0.75\linewidth]{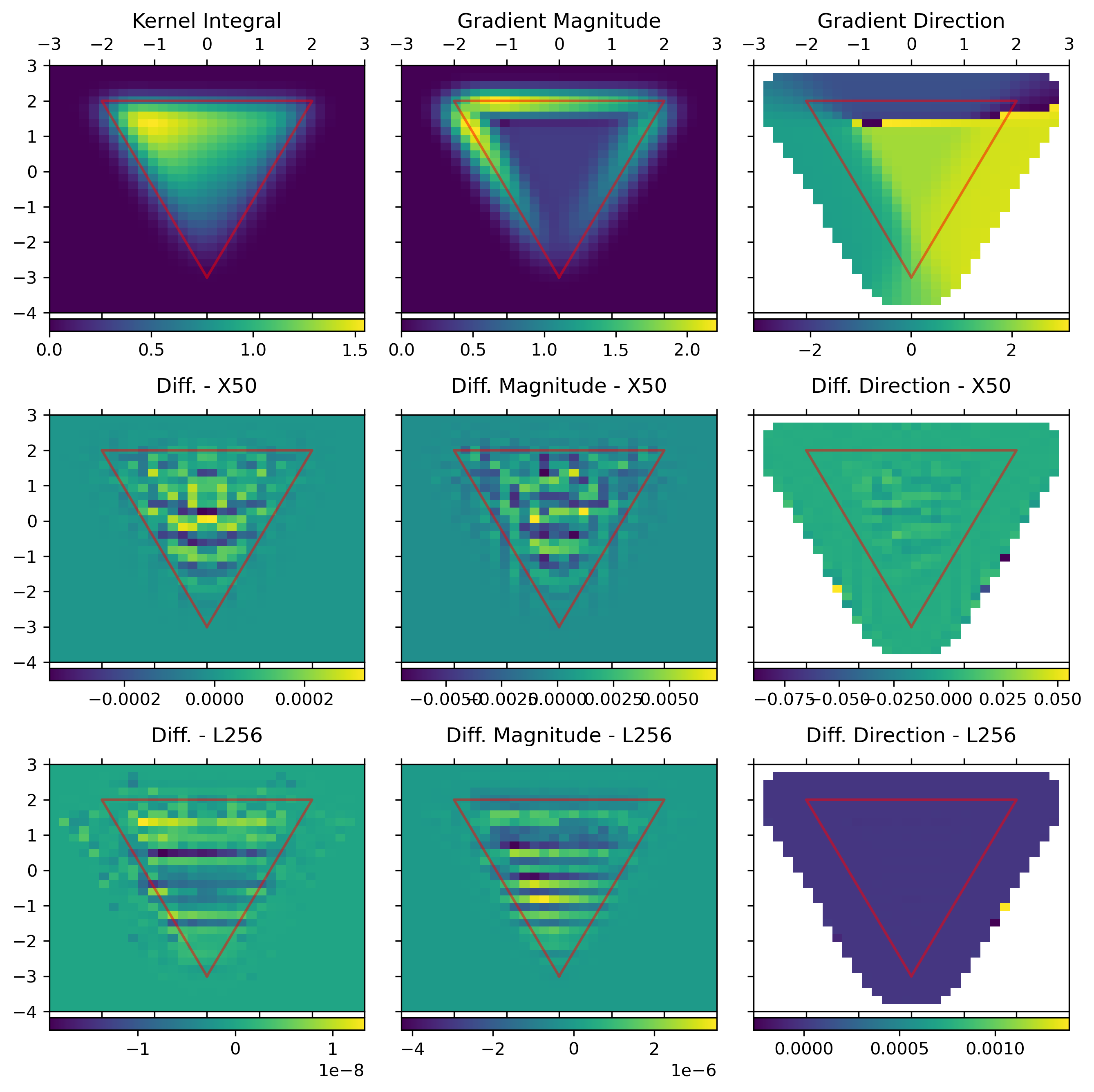}
    \caption{This figure shows the evaluation of the boundary (from left to right) integral as well as the magnitude and direction of the gradient for a triangle using piecewise-linear boundary quantities and a support $h=1$ for the integral. The top row shows the results evaluated using our analytic solution, the middle row shows the difference between our solution and using the $50$th order quadrature scheme of Xiao and Gimbutas~\cite{XIAO2010663} and the $256$th order scheme of Lether~\cite{LETHER1976219} using 65536 quadrature points.}
    \label{fig:largeTriangleExample}
\end{figure}

To qualitatively evaluate our proposed solutions against numerical schemes, we first compute the boundary integral values for the quantity itself and its gradient for a piecewise-constant field $f(x,y) = 1$ for a small triangle, i.e., a triangle where the longest edge is shorter than the support radius.
For such triangles, all integral domains will only partially cover the triangle, and accordingly, no trivial cases are present.
We highlight this case in Fig.~\ref{fig:smallTriangleExample}, which showcases good and smooth behavior of the overall solutions with values as expected from this case.
To actually evaluate the cases where numerical quadrature schemes have limits, we move to a slightly different test case.

In the second test case here, we evaluate the boundary integral solutions for a large triangle, i.e., one where the support radius is much smaller than the triangle, see Fig.~\ref{fig:largeTriangleExample}.
For these geometries, all possible cases discussed in Sec.~\ref{sec:integralTerms} occur throughout the evaluation domain.
Numerical quadratures for triangles cover the entire area of the triangle and utilize a polynomial approximation of the underlying field.
While the underlying field we consider in this case is linear, within the triangle, and the kernel function is a polynomial, the compact nature of the support domain of an SPH evaluation means that the function over the triangle, in this case, is piecewise-polynomial, and approximating this quantity is significantly harder using quadratures.
Furthermore, for cases where there is only partial overlap between the support domain and the triangle there might be few, or even no, quadrature points within the support domain.

Consequently, we expect the errors to focus on regions with only marginal overlap with the triangle due to the sparse nature of quadratures and to show the geometric distribution of the quadrature points when the full support domain is within the triangle.
Considering the results in Fig.~\ref{fig:largeTriangleExample}, we see both of these effects clearly as the areas of biggest error, with respect to the gradient direction, are for points that are only slightly closer than $h$ to the triangle, with a clear geometric pattern to the integral within the triangle.
Note that for small triangles that are fully within the support domain of an evaluation point, numerical quadrature is optimal, assuming the degree of the quadrature rule is appropriate, whereas for partial overlap, the compact nature of the SPH interpolant limits the achievable accuracy.
As a result, reducing the size of the triangles used to mesh a boundary geometry should lead to a convergence to the analytic solution as with smaller triangles a larger, relative, part of the domain is evaluated appropriately using quadrature with a smaller, relative, region being subject to the compact cutoff nature of SPH.
Consequently, a natural extension of our approach would be to utilize numerical quadrature if the triangle is fully within the support domain and use our analytic solution for all other cases; however, we will focus only on \emph{pure} solutions here.

\subsection{Quantitative Evaluation}
\begin{figure}
\centering
\begin{subfigure}{.45\textwidth}
  \centering
  \includegraphics[width=\linewidth]{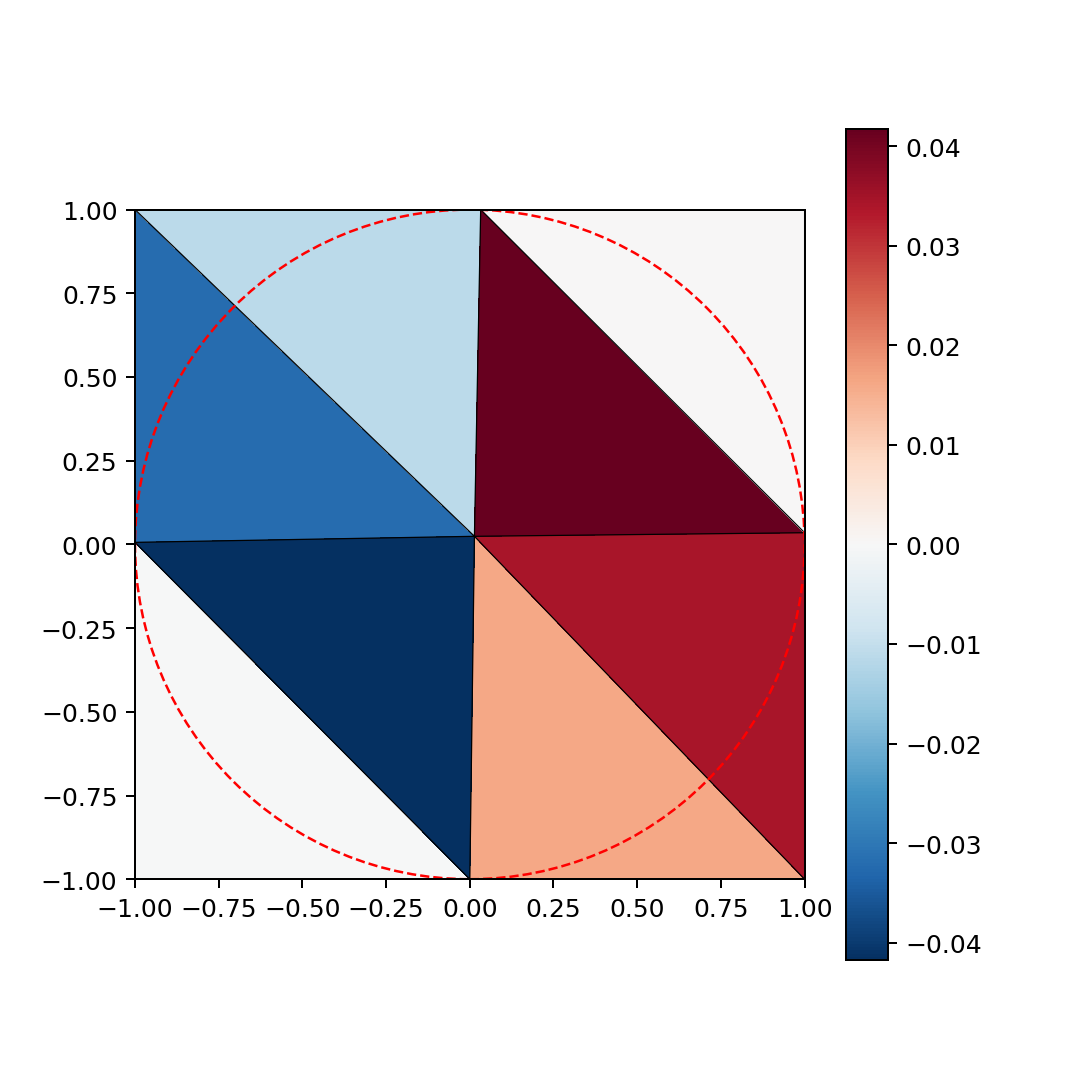}
  \caption{Result for $n_t = 8$.}
  \label{fig:eval:large_mesh}
\end{subfigure}
\begin{subfigure}{.45\textwidth}
  \centering
  \includegraphics[width=\linewidth]{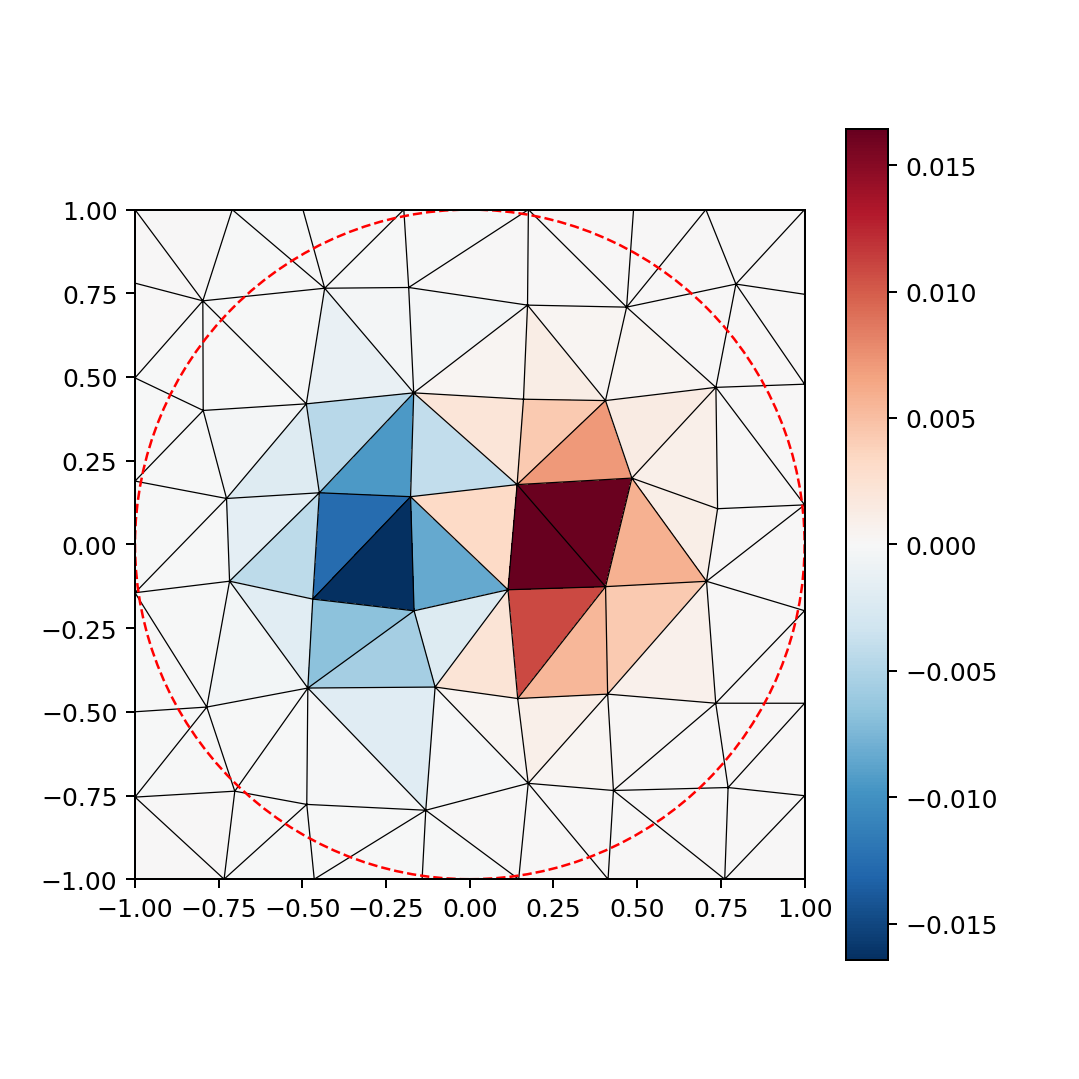}
  \caption{Result for $n_t = 96$.}
  \label{fig:eval:small_mesh}
\end{subfigure}%\\
\caption{Visualization of the result of evaluating the boundary integral over a set of triangles covering the $[-1,1]^2$ domain using a test function of $f(x,y) = x$ for two different numbers of triangles. Note that the visualized quantity is after convolving $f(x,y)$ with the kernel function.}
\label{fig:integralSmallLargeMesh}
\end{figure}

\begin{figure}
    \centering
    \includegraphics[width=1.0\linewidth]{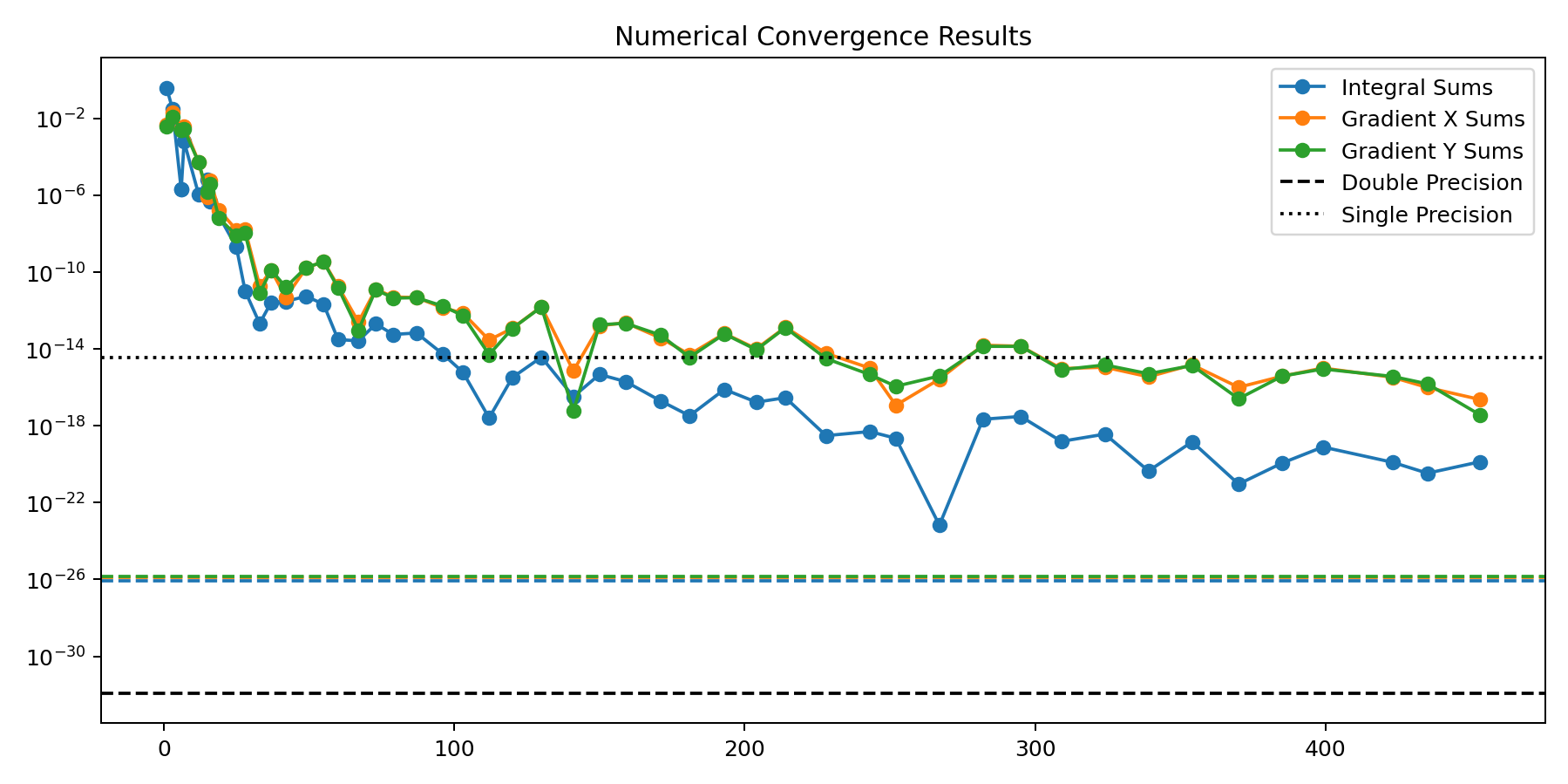}
    \caption{Convergence of the $L^2_2$ error of the numerical quadrature rule of Xiao and Gimbutas with increasing number of quadrature points per triangle compared to our analytic solution for the integral and gradients for the large triangle problem (see Fig.~\ref{fig:eval:small_mesh}) with piecewise-constant boundary elements. The black lines indicate machine precision epsilon values.}
    \label{fig:convergence_piecewise_constant}
\end{figure}
As discussed before, with numerical quadratures we expect the error of the integral to converge based on the size of the triangles, where we chose two test cases to evaluate this convergence and as test cases to evaluate our analytic solution as these cases have known solutions.
The first test case involves coupling a piecewise constant field with $f(x,y)=1$ over a coarse triangulation of the $[-1,1]^2$ domain with $8$ mesh elements in total and evaluating the boundary integral for $\mathbf{x}=\mathbf{0}$ with $h=1$, see Fig.~\ref{fig:convergence_piecewise_constant}, which should yield a boundary integral value of $1$ and a gradient of $\mathbf{0}$.
As the numerical quadrature scheme here we chose the quadrature rules of Xiao and Gimbutas~\cite{XIAO2010663}, as they allow us to generate quadrature rules of degree $2$ to $50$ with increasing quadruature points showing a convergence with respect to the order of quadrature.
In our evaluations, we observe a convergence of the numerical scheme towards the known solution at a slow convergence, yielding an $L_2$ error of $10^{-10}$ for the integral and $10^{-8}$, whereas our analytic solution yields $L_2$ errors of $10^{-13}$ for both.
While our analytic solution does not yield a solution that is exact within a machine precision epsilon, it is significant closer to machine precision than any numerical quadrature and the absolute numerical performance could be further improved by changes to the implementation.

\begin{figure}
    \centering
    \begin{subfigure}{.45\textwidth}
      \centering
      \includegraphics[width=\linewidth]{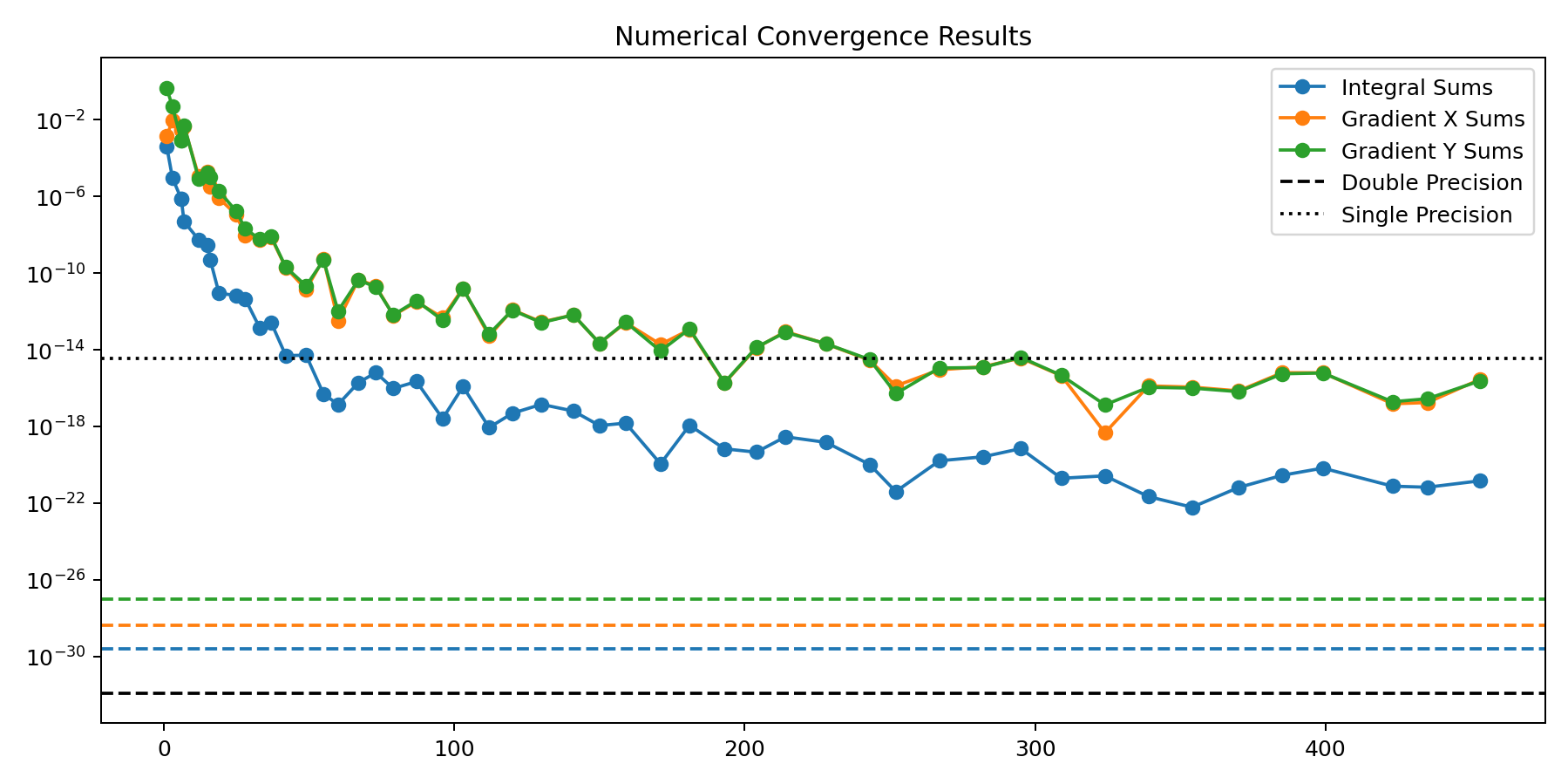}
      \caption{For large triangle elements, see Fig.~\ref{fig:eval:large_mesh}}
      \label{fig:convergence_linear_large}
    \end{subfigure}    
    \begin{subfigure}{.45\textwidth}
      \centering
      \includegraphics[width=\linewidth]{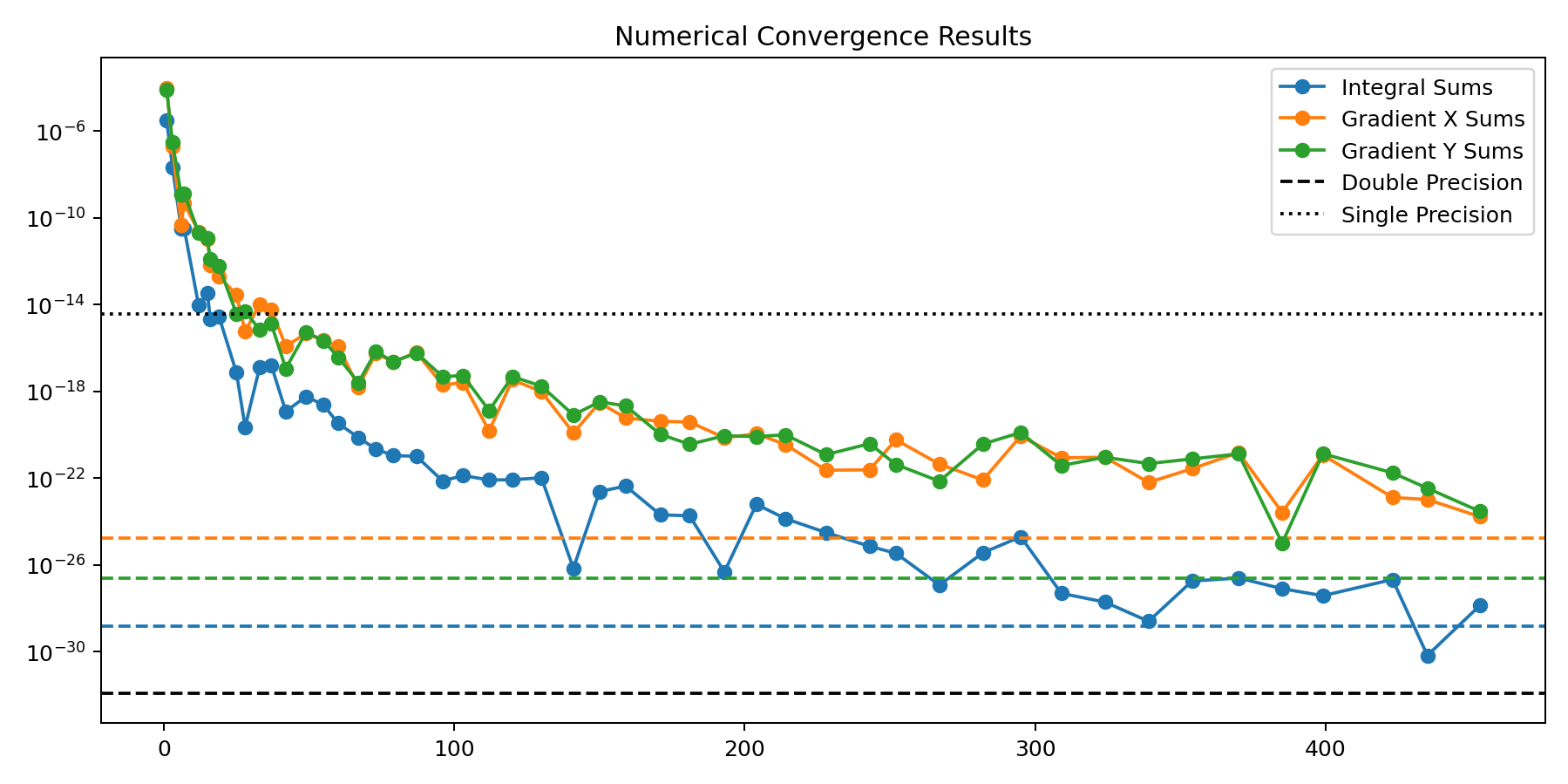}
      \caption{For small triangle elements, see Fig.~\ref{fig:eval:small_mesh}}
      \label{fig:convergence_linear_small}
    \end{subfigure}%
    \caption{Convergence of the $L^2_2$ error of the numerical quadrature rule of Xiao and Gimbutas with increasing number of quadrature points per triangle compared to our analytic solution for the integral and gradients for the small triangle problem (see Fig.~\ref{fig:eval:small_mesh}) with piecewise-linear boundary elements. The black lines indicate machine precision epsilon values.}
    \label{fig:convergence_piecewise_linear}
\end{figure}

As the second test, we then consider $f(x,y)=x$ with a boundary integral of $0$ and a gradient of $[1,0]^T$, for which we evaluate the errors both for a set of $8$ large triangles, see Fig.~\ref{fig:convergence_linear_large}, and a set of $96$ smaller triangles, see Fig.~\ref{fig:convergence_linear_small}.
For the large triangles we observe a very similar behavior for the numerical quadratures, as the underlying boundary element is only increase by $1$ order and the main difficulty being the compact nature of the SPH kernel, and our analytic solution remains largely unaffected.
For the smaller triangles we can see a significant improvement of the $L_2$ error for the numerical quadratures by $4$ orders of magnitude as the central region of the SPH kernel, with the largest values, is now represented by triangles that are fully within the support domain, whereas only a small fraction of triangles with small contributions is affected by the compact cutoff, see Fig.~\ref{fig:integralSmallLargeMesh}.
For our analytic solution we observe a slight decrease in the numerical accuracy, especially for the gradient terms, as significantly more numerical operations need to be performed here in total.
Overall our analytic solution still remains reliably more accurate than a numerical quadrature rule that required the evaluation of up to $3\times96\times 500 = 144000$ kernel functions to evaluate the quadrature points, whereas our analytic solution provided higher accuracy with only $8$ mesh elements.

\subsection{Incompressible SPH}
\begin{figure}
\centering
\begin{subfigure}{.33\textwidth}
  \centering
  \includegraphics[width=\linewidth]{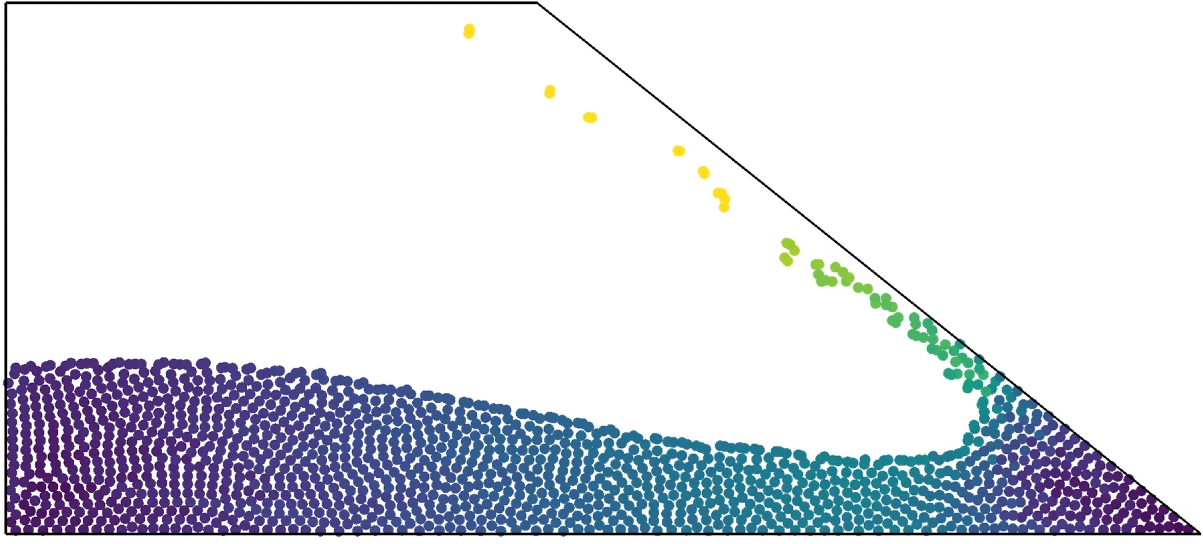}
  \caption{Acute corner}
  \label{fig:example_acute}
\end{subfigure}%
\begin{subfigure}{.33\textwidth}
  \centering
  \includegraphics[width=\linewidth]{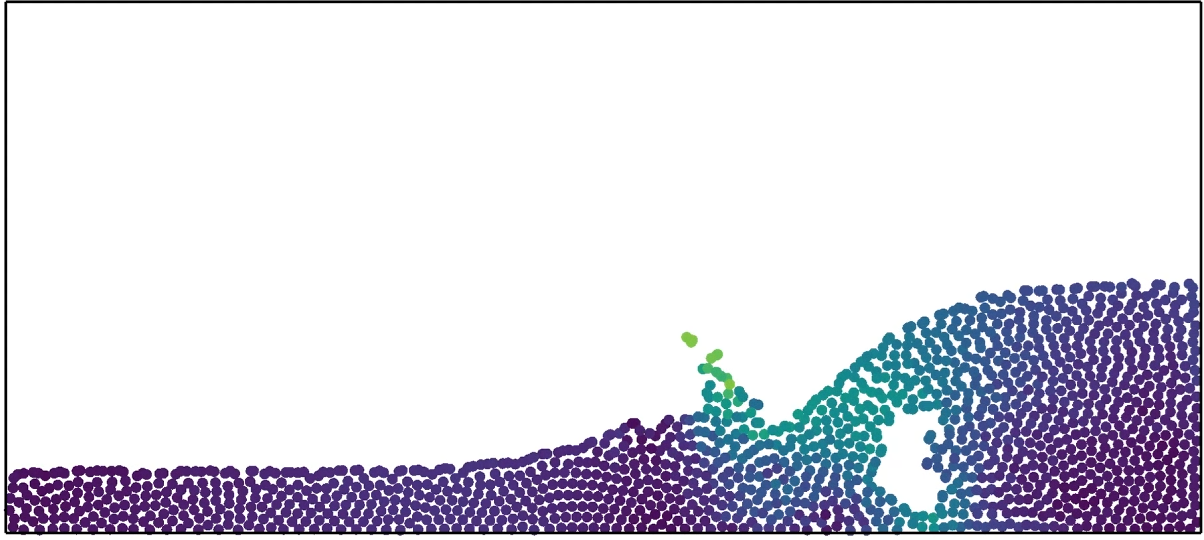}
  \caption{Orthogonal corner}
  \label{fig:example_ortho}
\end{subfigure}
\begin{subfigure}{.33\textwidth}
  \centering
  \includegraphics[width=\linewidth]{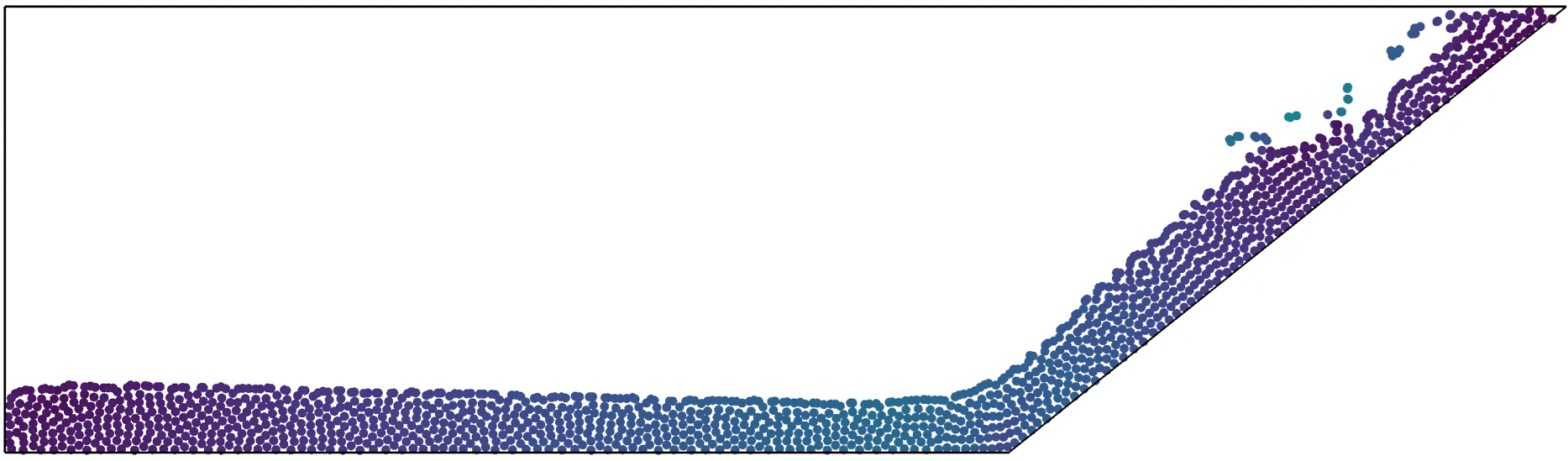}
  \caption{Obtuse corner}
  \label{fig:example_obtuse}
\end{subfigure}\\
\begin{subfigure}{.33\textwidth}
  \centering
  \includegraphics[width=\linewidth]{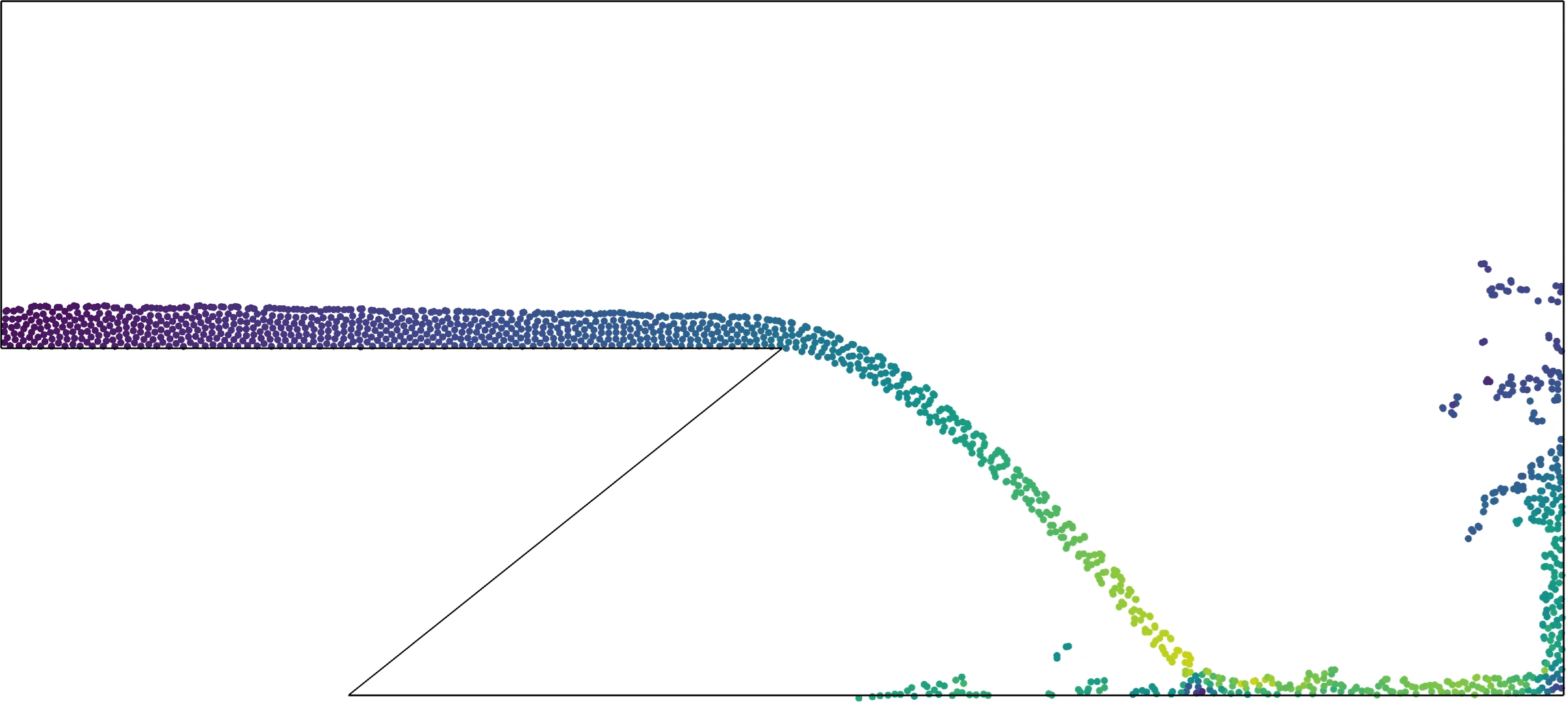}
  \caption{Acute corner with drop}
  \label{fig:example_nacute}
\end{subfigure}%
\begin{subfigure}{.33\textwidth}
  \centering
  \includegraphics[width=\linewidth]{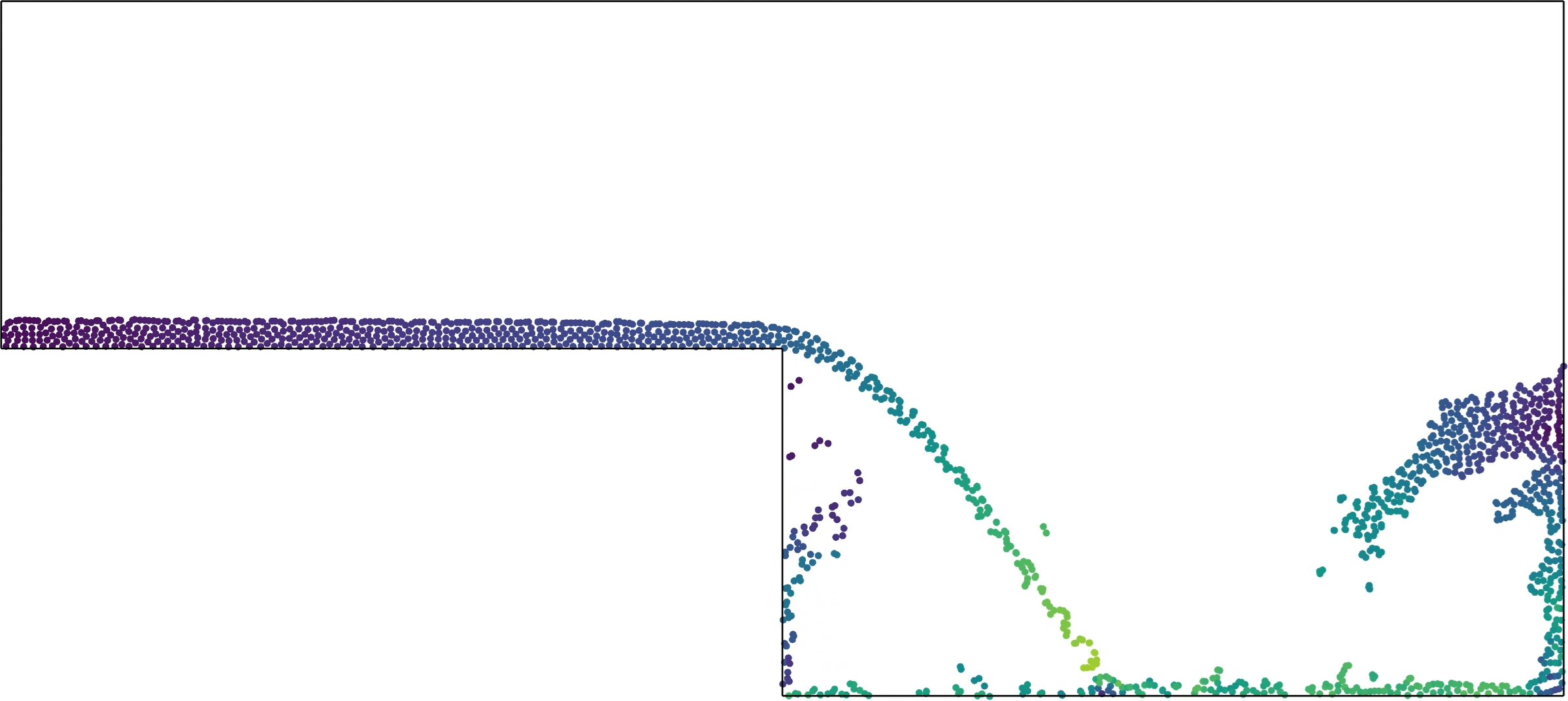}
  \caption{Orthogonal corner with drop}
  \label{fig:example_northo}
\end{subfigure}
\begin{subfigure}{.33\textwidth}
  \centering
  \includegraphics[width=\linewidth]{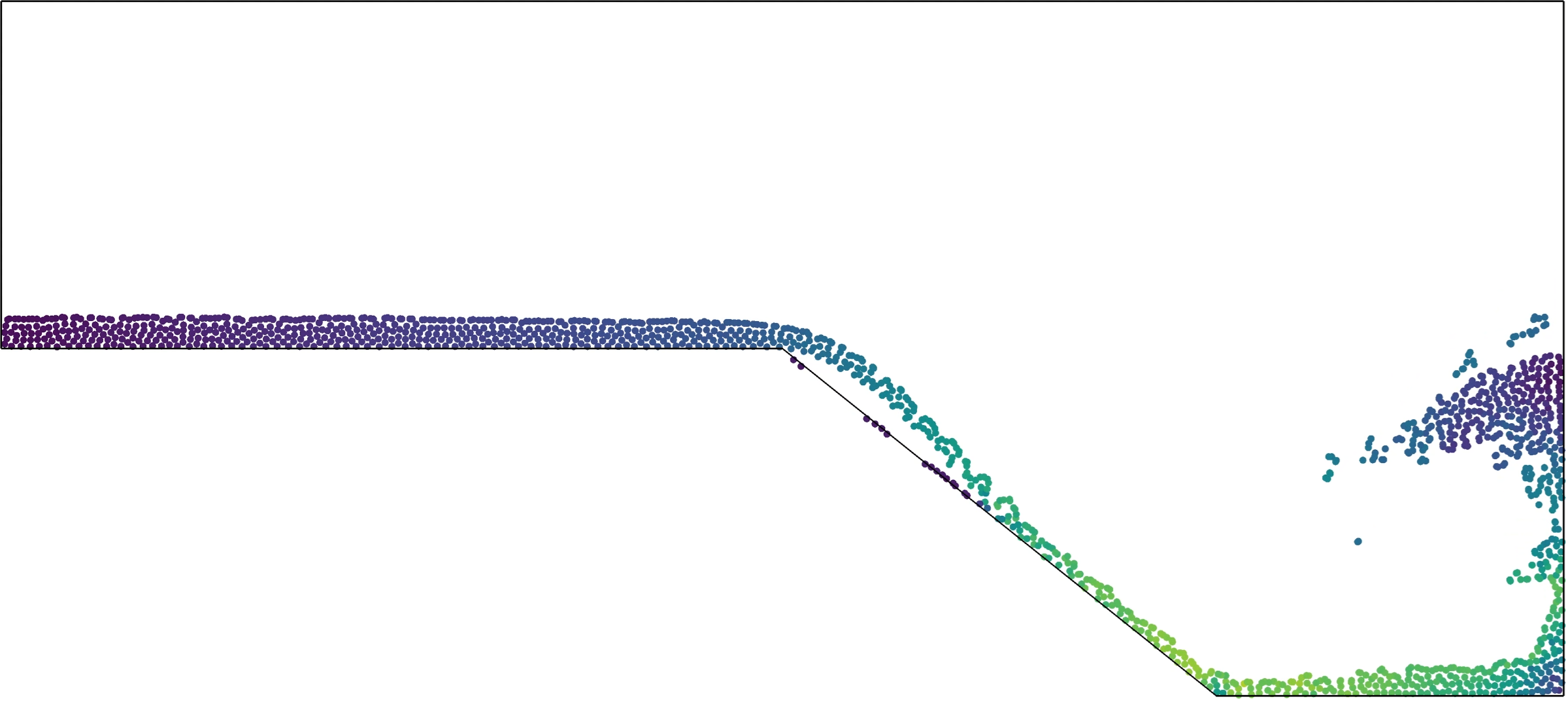}
  \caption{Obtuse corner with drop}
  \label{fig:example_nobtuse}
\end{subfigure}
\caption{Example simulations of using our boundary integral solution for an incompressible SPH simulation showing smooth flow behavior without artifacts in corners. Velocity color-coded.}
\label{fig:example_simulations}
\end{figure}

As our final test problem, we implement our boundary integral scheme in an incompressible SPH simulation, see Appendix~\ref{sec:SPH} for more details on the integration.
While developing fluid structure schemes that yield exact and robust results, especially in challenging conditions, is an important aspect of SPH, and more general CFD, research, our boundary integral solution is a general framework that can be used to develop these schemes.
Implementing our framework into such a broad range of schemes, including FEM coupling and beyond, is an important area of future research but beyond our scope here, as we focus on the solution of the integrals and how they can be directly integrated into an SPH scheme. \footnote{Note that an implementation of our SPH solver will be available publicly.}
This being said, for our test setup we utilize a set of six different boundary geometries that represent a dambreak interacting with a corner that either has an acute, orthogonal or obtuse angle, with half of the test cases including a fluid impact due to the step in the geometry.

In all of our test cases, see Fig.~\ref{fig:example_simulations}, we find a stable and robust simulation that does not yield any boundary penetration or artifacting.
These results are in close agreement with those of~\cite{DBLP:journals/tog/WinchenbachA020}, which utilize a numerical quadrature of $50$-th order to compare their semi-analytic solutions against.
When evaluating the pressure for the boundary elements we utilized two approaches as (a) computing the pressure only for nodes of the triangle mesh, and using piecewise-linear interpolation and (b) using a local contact point and evaluating the pressure for each contact point (note that this results in a different pressure for the same boundary element for each particle).

For the first approach we found that we needed to utilize very small boundary elements, on the order of $\frac{1}{4}h$, as the pressure field, especially for single particle impacts, is highly non-linear and requires either a small boundary element size, to accurately represent the non-linear field using piecewise-linear functions, or a high order boundary element.
This highlights the importance of choosing a correct discretization of the boundary element, to represent the required fluid mechanics, even if the evaluating of the coupling is exact.
For the second approach we found no dependence on the boundary element size and were able to utilize very few triangles for the entire domain, e.g., the acute corner impact case in Fig.~\ref{fig:example_acute}, consisted of a total of $8$ triangles.
Our analytic boundary integral in this case yielded exact gradient and kernel interpolation magnitudes, while the local contact point formulation yielded an appropriate pressure value for each particle.
Overall this clearly highlights the potential of our analytic boundary integrals and highlights several key directions for future research.

\section{Conclusion}
\label{sec:conclusion}
In this paper, we presented an analytic and algebraic solution to couple arbitrary (piecewise) order boundary quantities defined on a triangular mesh exactly to an SPH simulation for arbitrary polynomial kernels.
Using our analytic solutions, we achieve consistent results across varying resolutions and regardless of triangle geometry, which outperform numerical quadrature-based solutions by up to 5 orders of magnitude, especially for larger triangle elements.
As our solution directly models the volume contribution of the boundary, instead of relying on surface integrals, the solution can be applied not solely to gradient quantities, but also directly for SPH interpolants, enabling a seamless and direct integration of SPH with FEM.
While using SPH interpolants made it possible to exactly couple from SPH to FEM, our solution closes the gap and enables exact and analytic two-way coupling between FEM and SPH, enabling a promising future avenue of research and coupling.

While our method is analytic and exact, it is not without drawbacks.
The evaluation of some of the integral terms can be computationally challenging, especially for GPUs, due to the large number of branching operations.
Furthermore, our solution is only applicable to two-dimensional problems, and directly applying our solution to tetrahedra in three dimensions is a significant but promising challenge.
Finally, studying the relationship between order of coupling, size of triangle elements, and the resulting fluid behavior is beyond the scope of our work here, but our framework builds the necessary foundations to perform these studies in the future.
Note that an implementation of our framework, both in Python and C++, will be made available publicly.

\appendix
\section{Integration into Pressure Solvers}

To enforce incompressibility, we follow the approach of Bender and Koschier~\cite{bender2015divergence}, based on Ihmsen et al.~\cite{ihmsen2013implicit}, where a pressure projection method is used to first ensure that the velocity field is divergence free and then to ensure that the velocity field is incompressible.
In this approach the continuity equation $\frac{D\rho}{Dt}=-\rho\nabla\cdot v$ is discretized using a forward difference for the material derivative of density, i.e., $\frac{D\rho}{Dt} = \frac{\rho^{t+\Delta t}-\rho^t}{\Delta t}$ and a difference formulation for the divergence term $\nabla\cdot v$, yielding
\begin{equation}
\begin{split}
\frac{\rho_i(t+\Delta t)-\rho}{\Delta t} = &\sum_j m_j \mathbf{v}_{ij}^\text{adv}\cdot\nabla W_{ij} \\ &+ \Delta t \sum_j m_j \left[\frac{\mathbf{F}_i^p}{m_i} - \frac{\mathbf{F}_j^p}{m_j}\right]\cdot \nabla W_{ij},
\end{split}
\end{equation}
with $\mathbf{v}_{ij} = \mathbf{v}_{i} - \mathbf{v}_{j}$, $\nabla W_{ij} = \nabla_{\mathbf{x}_i}W(\mathbf{x}_i - \mathbf{x}_j,h)$ and $\mathbf{F}_i^p = -\sum_j m_j \left[\frac{p_i}{\rho_i^2} + \frac{p_j}{\rho_j^2}\right]\nabla W_{ij}$. This yields a system of equations with unknown pressure values $p$ for all particles and an unknown left hand side of the equation system. To solve this system for incompressibility the density in the new timestep is assumed to be equal to the rest density~\cite{ihmsen2013implicit,band2018pressure}, i.e., $\rho_i(t+\Delta t) = \rho_0$, and to solve for divergence-freedom the left hand side is assumed to be equal to 0~\cite{band2018mls,bender2015divergence}.

To evaluate the acceleration of a particle $i$ due to pressure forces from other particles, $\mathbf{a}_i^f$, it is necessary to evaluate the pressure gradient as~\cite{monaghan1994simulating}
\begin{equation}
    \mathbf{a}_i^f = - \frac{\nabla p_i}{\rho_i} \approx -\sum_j m_j \left(\frac{p_i}{\rho_i^2} + \frac{p_j}{\rho_j^2}\right)\nabla_i W_{ij},
\end{equation}
where a symmetric gradient formulation was used to ensure symmetric pressure forces.
Evaluating the acceleration of a particle due to pressure forces from a boundary region $\Omega$, $\mathbf{a}_i^b$, involves evaluating the integral~\cref{eqn:intsphdiffops} by inserting the pressure field $p(x)$ for $A$, which yields
\begin{equation}
   \mathbf{a}_i^b  = -\frac{\nabla p}{\rho} \approx -\int_\Omega  \rho(\mathbf{x}^\prime)\left[\frac{p(\mathbf{x}^\prime)}{\rho(\mathbf{x}^\prime)^2}+\frac{p_i}{\rho_i^2}\right]\nabla_i W(\mathbf{x} - \mathbf{x}^\prime,h)d\mathbf{x}^\prime.
\end{equation}

Where, based on the prior derivation, see Sec.~\ref{sec:integralTerms}, we can evaluate this term for a scalar field $f(\mathbf{x}^\prime) = \rho(\mathbf{x}^\prime)\left[\frac{p(\mathbf{x}^\prime)}{\rho(\mathbf{x}^\prime)^2}+\frac{p_i}{\rho_i^2}\right]$ as 
\begin{equation}
    \mathbf{a}_i^b \approx -\int_\Omega f(\mathbf{x}^\prime) \nabla_i W(\mathbf{x} - \mathbf{x}^\prime,h)d\mathbf{x}^\prime,
\end{equation}
which yields a total acceleration due to pressure forces of a particle $i$ due to neighboring particles $j$ and neighboring triangles $\mathcal{T}_j$ as
\begin{equation}
\begin{split}
\mathbf{a}_i = &-\sum_j m_j \left(\frac{p_i}{\rho_i^2} + \frac{p_j}{\rho_j^2}\right)\nabla_i W_{ij}\\
&-\sum_{\mathcal{T}_j}\int_{\mathcal{T}_j} f(\mathbf{x}^\prime) \nabla_i W(\mathbf{x} - \mathbf{x}^\prime,h)d\mathbf{x}^\prime.
\end{split}
\end{equation}

In our implementation we assume that the density of the boundary domain is equal to the rest density of the boundary object, i.e., $\rho(\mathbf{x}^\prime) = \rho_{0,b}$, as we assume boundary objects to be non-deformable and, consequently, their density cannot deviate from their rest density.

Regarding the pressure on the boundary object, $p(\mathbf{x}^\prime)$, we propose two distinct choices to model this term. 
First, a pressure value at each vertex of a triangle could be determined and then barycentrically interpolated across the triangle region and second, a pressure value for the entire triangle could be determined with a distinct value per fluid particle.
The first choice involves evaluating a pressure value for each vertex, i.e., $p_1 = p({\nu_1}), p_2 = p(\nu_2)$ and $p_3 = p({\nu_3})$ for each triangle, whereas the second choice involves evaluating a pressure value per particle for each triangle $p(\mathbf{x}_{i,\mathcal{T}})$, where $\mathbf{x}_{i,\mathcal{T}}$ is the closest point on the boundary surface $\delta\mathcal{T}$ relative to the position of the particle $\mathbf{x}_i$.

To evaluate a pressure value at an arbitrary position $\mathbf{x}$ we utilize the approach of Band~et. al~\cite{band2018mls}, where the pressure value is evaluated as the result of a Moving Least Squares interpolation.
Consequently, for the first choice we need to evaluate a barycentrically interpolated term using the three pressure values $p_1, p_2$ and $p_3$ as
\begin{equation}
    \mathbf{a}_i^b \approx \iint_\mathcal{T} 
     \left[\begin{array}{c}\tau_1^f g_x^1(r,\theta) + \tau_2^f g_x^2(r,\theta) + \tau_3^f g_x^3(r,\theta)\\\tau_1^f g_y^1(r,\theta)+\tau_2^f g_y^2(r,\theta)+\tau_3^f g_y^3(r,\theta)\end{array}\right]d\theta dr,
\end{equation}
with $\tau^f$ defined as before, see Sec.~\ref{sec:integrals}, whereas for the second choice the field is constant across the triangle and, consequently, the term can be simplified to yield
\begin{equation}
    \mathbf{a}_i^b \approx \rho_{0,b}\left[\frac{p(\mathbf{x}_{i,\mathcal{T}})}{\rho_{0,b}^2} + \frac{p_i}{\rho_i^2}\right]\iint_\mathcal{T} 
    \left[\begin{array}{c}g_x^3(r,\theta)\\g_y^3(r,\theta)\end{array}\right]d\theta dr.
\end{equation}

Accordingly, the first choice requires an evaluation of three pressure values, using MLS, per triangle and a barycentric approach per particle, whereas the second choice requires an evaluation of one pressure value, using MLS, per neighboring triangle per particle and a non barycentric approach per particle.
These terms can then be directly integrated into DFSPH~\cite{DBLP:journals/tog/WinchenbachA020}.

To model boundary terms we follow an approach similar to that of Winchenbach~et al.~\cite{DBLP:journals/tog/WinchenbachA020}, where boundary conditions are modelled by assuming a local planar boundary region.
Evaluating boundary condition terms on a per triangle mesh element basis would not yield desirable results as the boundary normals of individual mesh elements does not reflect the overall boundary normal resulting in spurious behavior. 
To resolve this issue we first evaluate a boundary normal, i.e., we evaluate
\begin{equation}
    \mathbf{n}_i = \int_\mathcal{T} \nabla W(\mathbf{x}_i - \mathbf{x}^\prime,h)d\mathbf{x}^\prime,
\end{equation}
using our proposed analytic boundary integration scheme.
It is important to note that this has to be evaluated for any boundary object region, e.g.,  a particle flowing through a narrow gap inside of a boundary object would otherwise yield an incorrect boundary normal of $\mathbf{0}$.
We then evaluate the magnitude of the viscosity term by integrating over all triangle mesh elements and apply only the tangential part, relative to the boundary normal evaluate above, to the particle velocity.
It is important to note that most boundary conditions, e.g., a standard Navier-Stokes viscosity term, given by
\begin{equation}
    \frac{d\mathbf{v_i}}{d t} = m_i\nu \nabla^2 \mathbf{v}_i,
\end{equation}
cannot be evaluated in a barycentric boundary formulation.
In general, evaluating a direct second order derivative term leads to a force that is not momentum conserving~\cite{price2012smoothed} and instead a finite difference scheme is applied to a first order derivative~\cite{brookshaw1985method,price2012smoothed,SPHTutorial}, which yields a discretized laplacian as
\begin{equation}
    \nabla^2\mathbf{A}_i \approx 2(d+2)\sum_{j\in\mathcal{N}_i} \frac{m_j}{\rho_j} \frac{\mathbf{A}_{ij}\cdot\mathbf{x}_{ij}}{||\mathbf{x}_{ij}||^2}\nabla_i W_{ij}.
\end{equation}
Applying this formulation to the velocity field then yields a viscosity term~\cite{monaghan2005smoothed}
\begin{equation}
    \frac{d\mathbf{v_i}}{d t} = m_i\nu 2(d+2)\sum_{j\in\mathcal{N}_i} \frac{m_j}{\rho_j} \frac{\mathbf{v}_{ij}\cdot\mathbf{x}_{ij}}{||\mathbf{x}_{ij}||^2 + 0.01 h^2}\nabla_i W_{ij},
\end{equation}
where the term $0.01 h^2$ was added to the denominator to avoid singularities.
It is possible to treat this quantity as being linearily interpolated across a boundary element, however, this term is not linear in nature and would yield nonphysical boundary interactions.
Accordingly, we opt to only interact with a single planar representation of the boundary, for velocity boundary conditions, with a single contact point~\cite{DBLP:journals/tog/WinchenbachA020}.

\section*{Acknowledgements}
Funding was provided by the German Research Foundation~(DFG) TH 2034/1-2.


\bibliographystyle{unsrt}

\bibliography{sn-bibliography}



\end{document}